\documentclass[11pt]{amsart}
\usepackage{fullpage}
\usepackage{amsmath,amssymb,amsthm}
\usepackage[foot]{amsaddr}
\usepackage[colorlinks=true,linkcolor=red,citecolor=blue,urlcolor=cyan]{hyperref}
\usepackage{xcolor}

\newtheorem{theorem}{Theorem}
\newtheorem{proposition}{Proposition}

\usepackage{graphicx}
\usepackage{caption}
\usepackage{booktabs}
\usepackage{siunitx}

\usepackage[backend=biber,giveninits=true,sorting=nyt,style=alphabetic,natbib=true,maxcitenames=2,maxbibnames=10,url=false,doi=true,backref=false]{biblatex}
\renewbibmacro{in:}{\ifentrytype{article}{}{\printtext{\bibstring{in}\intitlepunct}}}

\allowdisplaybreaks
\AtBeginDocument{\setlength\abovedisplayskip{4pt plus 1pt minus 1pt}\setlength\belowdisplayskip{4pt plus 1pt minus 1pt}\setlength\abovedisplayshortskip{1pt}\setlength\belowdisplayshortskip{2pt}}
\begin{document}

\author{Haoyu Chen}
\address{Department of Mathematics, University of Utah, Salt Lake City, UT}
\email{Haoyu.Chen@utah.edu}
\thanks{The authors acknowledge partial support from NSF DMS-2136198 and DMS-2513175.}

\author{Braxton Osting}
\email{osting@math.utah.edu}

\title{Computation of free boundary minimal \\ tori via extremal Steklov eigenvalues}



\date{\today}

\begin{abstract}
Fraser and Schoen showed that the first Steklov eigenfunctions of a metric maximizing the length-normalized first eigenvalue of a surface with boundary define a free boundary minimal immersion into a Euclidean ball. We develop a method for this extremal problem on surfaces of genus one with \(b\) boundary components and use it to compute free boundary minimal tori for \(b\in\{1,2,3,4,5,8\}\). By uniformization, the problem is a nonsmooth maximization over the modulus of a flat torus, the centers and radii of \(b\) removed disks, and a boundary weight. The eigenvalues are computed with doubly periodic layer potentials, all derivatives follow from a single formula for the underlying pencil of boundary integral operators, and the coalescence of eigenvalues at a maximizer is handled by a two-stage nonsmooth method. The Lagrange multiplier matrix that certifies first-order stationarity rescales the eigenfunctions into the free boundary minimal immersion. For \(b=1\), the supremum is the Fraser--Schoen bound \(4\pi\), attained by branched double covers of the disk; the embedded torus with connected boundary of Carlotto, Franz, and Schulz is a saddle point, which answers their question whether it maximizes the eigenvalue in the negative. For \(b\geq2\) the largest values correspond to free boundary minimal tori in \(\mathbb B^5\) and \(\mathbb B^6\); the surfaces of Folha, Pacard, and Zolotareva, two disks joined by a central neck and \(b\) equatorial necks, are maximizers within their symmetry class but saddle points of the full problem, with values increasing toward \(4\pi\) as \(b\to\infty\).
\end{abstract}

\maketitle

\section{Introduction}
\label{sec:intro}

Extremal eigenvalue problems originate with Rayleigh's conjecture that, among planar domains of fixed area, the disk minimizes the principal Dirichlet--Laplace eigenvalue. In the century since its proof by Faber and Krahn, many generalizations have been studied \cite{Henrot_2006}, and for problems without an explicit solution the subject has become a challenging computational one \cite{Oudet_2004,Osting_2010,Antunes_2012,Osting_Kao_2013,Bogosel_Oudet_2016,Bogosel_2017,Akhmetgaliyev_2017}. Two features of the computed optimizers recur: they are highly symmetric, and the extremal eigenvalue is multiple \cite{Nadirashvili_1988}. Multiplicity makes the objective nondifferentiable exactly at the points sought, so stationarity must be understood in Clarke's sense and the optimization method must be built for it \cite{Overton_1988,Overton_Womersley_1993,Lewis_Overton_1996}, either by general-purpose nonsmooth methods such as gradient sampling and nonsmooth quasi-Newton iterations \cite{Burke_Lewis_Overton_2005,Lewis_Overton_2013,Curtis_Mitchell_Overton_2017} or, more efficiently, by methods that identify the active eigenvalue cluster and exploit the smooth structure of the problem on it \cite{Bareilles_2023}.

In an \emph{extremal geometric eigenvalue problem}, the metric rather than the domain is the variable: one extremizes a normalized eigenvalue of the Laplace--Beltrami or Dirichlet-to-Neumann operator over the Riemannian metrics on a fixed surface. The derivative of an eigenvalue with respect to the metric is classical \cite{ElSoufi_Ilias_2008}, and in a few cases the optimal metric is known and has a distinguished structure, such as the round sphere \cite{Hersch_1970} or the equilateral torus \cite{Nadirashvili_1996}. Where it has no such structure, computation is the available tool, and it faces two difficulties absent from the Euclidean setting: the space of conformal classes must be parametrized, with derivatives of the eigenvalue with respect to those parameters, and eigenvalue problems must be solved accurately on surfaces with nontrivial genus and boundary \cite{Kao_2022}. Computational work has so far treated closed surfaces of genus zero and one \cite{Kao_2017} and surfaces with boundary of genus zero \cite{Kao_2020}; extending it to surfaces with boundary and positive genus, where the modulus and the boundary weight enter together, was listed as open in \cite{Kao_2022}. This paper treats the first such case, the genus-one Steklov problem, whose maximizers are free boundary minimal surfaces.

Let \((\Sigma,g)\) be a smooth, compact, connected, orientable Riemannian surface of genus \(\gamma\) with \(b\geq1\) boundary components. The Steklov eigenproblem is to find eigenpairs \((\sigma,u)\) satisfying
\begin{equation}\label{eq:geometric_steklov_problem}
\begin{aligned}
\Delta_g u&=0&&\text{in }\Sigma,\\
\partial_{\nu_g}u&=\sigma u&&\text{on }\partial\Sigma.
\end{aligned}
\end{equation}
The eigenvalues, counted with multiplicity,  form a sequence \(0=\sigma_0<\sigma_1\leq\sigma_2\leq\cdots\to\infty\) and equivalently are  eigenvalues of the Dirichlet-to-Neumann operator \cite{Girouard_2017}. Since \(\sigma_k\) scales like an inverse length, the scale-invariant quantity is
\begin{equation}\label{eq:intro-normalized-eigenvalue}
\bar\sigma_k(\Sigma,g):=\sigma_k(\Sigma,g)\,L_g(\partial\Sigma),
\end{equation}
where \(L_g(\partial\Sigma)\) is the boundary length, and the extremal Steklov problem is
\begin{equation}\label{eq:intro-extremal-steklov-problem}
\sigma_k^*(\gamma,b):=\sup_g\ \bar\sigma_k(\Sigma,g),
\end{equation}
the supremum over all smooth metrics on \(\Sigma\); we write \(\sigma^*:=\sigma_1^*\). The problem goes back to Weinstock \cite{Weinstock_1954}, who proved \(\bar\sigma_1\leq2\pi\) for simply connected planar domains, with equality only for the disk. In general the supremum is finite, \(\bar\sigma_k\leq2\pi(\gamma+b+k-1)\) \cite{Fraser_2011,Girouard_2012,Karpukhin_2017}, and for \(k=1\) it is attained by a smooth metric \cite{Fraser_2015,Matthiesen_Petrides_2020}.

What makes the Steklov problem special among extremal geometric eigenvalue problems is a surprising connection between its maximizers and minimal surfaces, discovered by Fraser and Schoen \cite{Fraser_2011,Fraser_2015}. A \emph{free boundary minimal surface in the unit ball \(\mathbb B^m\)} is a critical point of area among surfaces with boundary on the sphere \(\mathbb S^{m-1}\), equivalently a minimal surface meeting \(\mathbb S^{m-1}\) orthogonally along its boundary. The equatorial disk and the critical catenoid are the classical examples, and constructing such surfaces with prescribed topology is challenging \cite{li2020}.
If \(X\colon\Sigma\to\mathbb B^m\) is a free boundary minimal immersion and \(g=X^*g_{\mathbb R^m}\), then the coordinate functions satisfy
\begin{equation}\label{eq:intro-coordinate-steklov}
\Delta_gX=0\ \ \text{in }\Sigma,\qquad
\partial_{\nu_g}X=X\ \ \text{on }\partial\Sigma,
\end{equation}
that is, they are Steklov eigenfunctions with eigenvalue one. Fraser and Schoen proved a converse: if \(g\) maximizes \(\bar\sigma_1\), then, after scaling, a suitable basis of the first eigenspace defines a possibly branched free boundary minimal immersion of \((\Sigma,g)\) into a ball, conformal with respect to \(g\), whose dimension is the multiplicity \(m\) of \(\sigma_1\); we call the map by first eigenfunctions the \emph{spectral map}. The correspondence has become a source of existence theorems for free boundary minimal surfaces \cite{Fraser_2015,Matthiesen_Petrides_2020,Karpukhin_Stern_2021,karpukhin2024}; read constructively, it turns eigenvalue optimization into a computational method, a maximizing metric with \(m=3\) yielding a free boundary minimal surface in \(\mathbb B^3\). The normalized eigenvalue has a geometric meaning on that surface: by \eqref{eq:intro-coordinate-steklov} and conformality, the boundary of the surface has length exactly \(\bar\sigma_1\) (Section~\ref{ssec:stationarity}). The extremal problem thus seeks the free boundary minimal surface of the given topology with the longest boundary, among those whose coordinates are first eigenfunctions, and the bound \(\bar\sigma_1\leq2\pi(\gamma+b)\) says that the boundary of a free boundary minimal torus with \(b\) boundary curves whose coordinates are first eigenfunctions is at most \(b+1\) times as long as a great circle; the equatorial disk and the critical catenoid have boundary length \(2\pi\) and \(10.47\).

The optimization problem solved in this paper is the genus-one case of \eqref{eq:intro-extremal-steklov-problem} with \(k=1\). By uniformization, a genus-one surface with \(b\) boundary components is conformally a flat torus \(\mathbb T_\tau=\mathbb C/(\mathbb Z+\tau\mathbb Z)\) with \(b\) round disks removed, and since \(\bar\sigma_1\) is invariant under conformal changes of the metric in the interior, the metric is described by a positive weight \(\rho\) on the boundary (Section~\ref{ssec:circle-domain}). With \(\log\rho\) truncated to \(N_\rho\) Fourier modes on each circle, the problem becomes
\begin{equation}\label{eq:intro-finite-dim-problem}
\max_{\mathbf p}\ \bar\sigma_1(\mathbf p),\qquad
\mathbf p=(\tau,\mathbf c,\mathbf r,\boldsymbol\alpha,\boldsymbol\beta)\in\mathbb R^{P},\qquad P:=2+3b+b(2N_\rho+1),
\end{equation}
over the modulus \(\tau\), the centers \(\mathbf c\) and radii \(\mathbf r\) of the disks and the Fourier coefficients \((\boldsymbol\alpha,\boldsymbol\beta)\) of \(\log\rho\), subject to the disks being pairwise disjoint. Problem \eqref{eq:intro-finite-dim-problem} is nonconvex, each evaluation of the objective requires solving a weighted Steklov eigenproblem on a finitely connected torus, and the objective is nonsmooth at the maximizers, where \(\sigma_1\) is expected to be multiple. The genus-zero analogue was treated by \textcite{Kao_2020}, who combined the conformal uniformization of planar domains with circular holes, a boundary-weighted Steklov eigensolver, and nonsmooth optimization to compute genus-zero free boundary minimal surfaces with many boundary components. The passage to genus one introduces the modulus as an optimization variable and requires an accurate eigensolver on the torus; we use the layer-potential theory for finitely connected tori of \textcite{kim2025layer}, which rests on the logarithmic conjugation theorem for harmonic functions on such domains \cite{KOO_Harmonic_2023}.

\begin{figure}[p]
    \centering
    \includegraphics[width=0.202\textwidth,trim={30 30 14 14},clip]{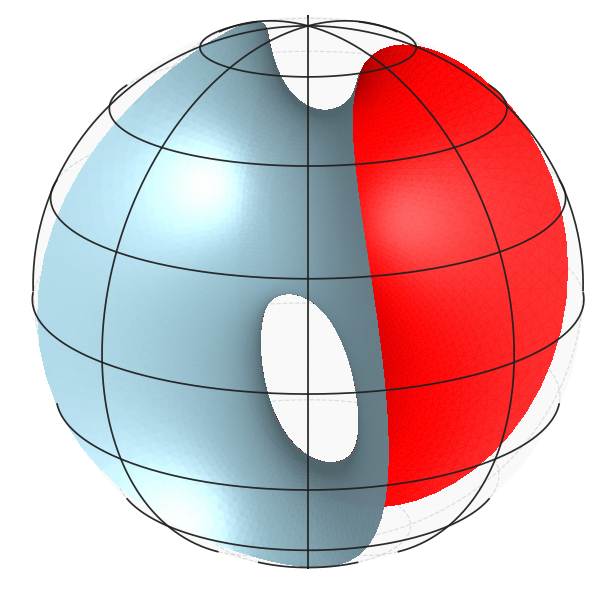}\hspace{2cm}
    \includegraphics[width=0.202\textwidth,trim={30 30 14 14},clip]{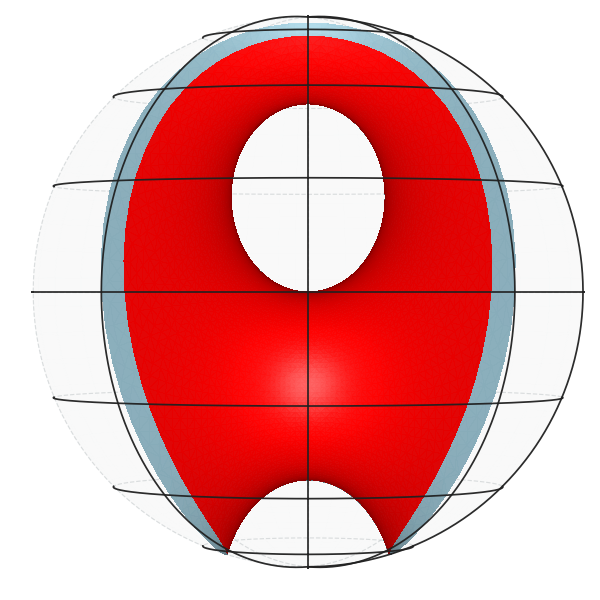}\hspace{2cm}
    \includegraphics[width=0.202\textwidth,trim={30 30 14 14},clip]{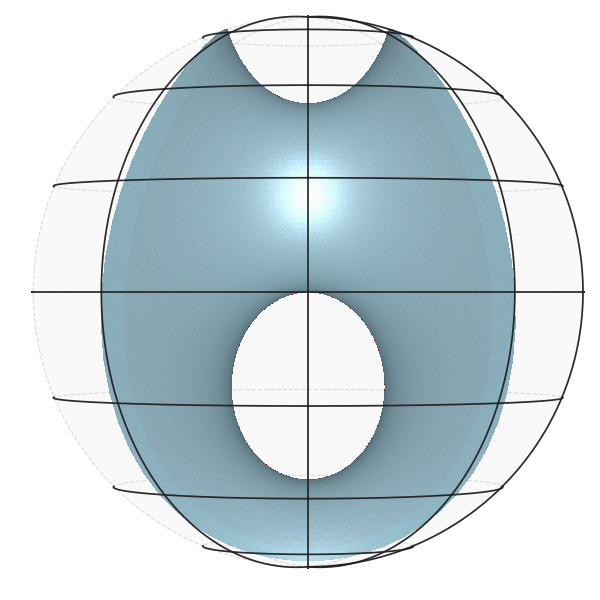}\\
    \includegraphics[width=0.202\textwidth,trim={30 30 14 14},clip]{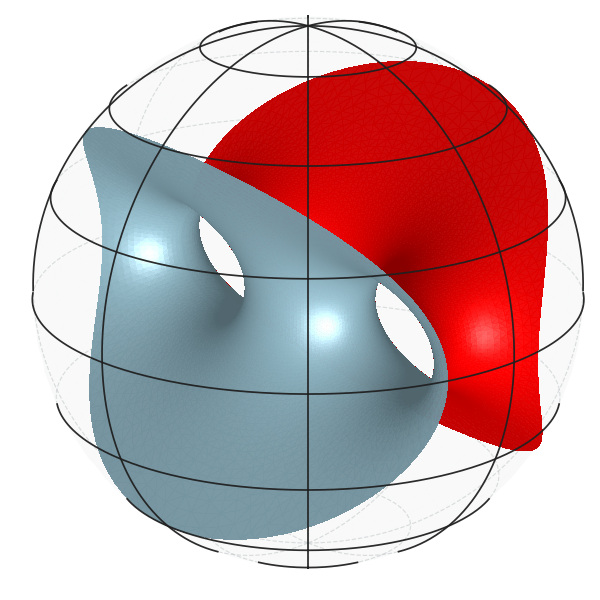}\hspace{2cm}
    \includegraphics[width=0.202\textwidth,trim={30 30 14 14},clip]{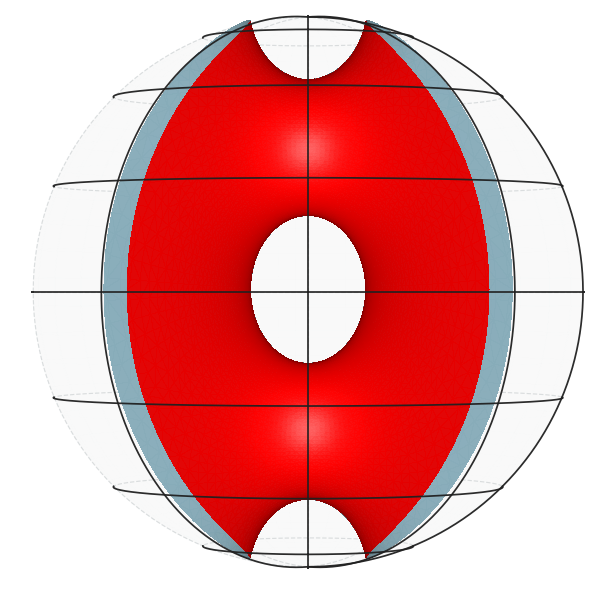}\hspace{2cm}
    \includegraphics[width=0.202\textwidth,trim={30 30 14 14},clip]{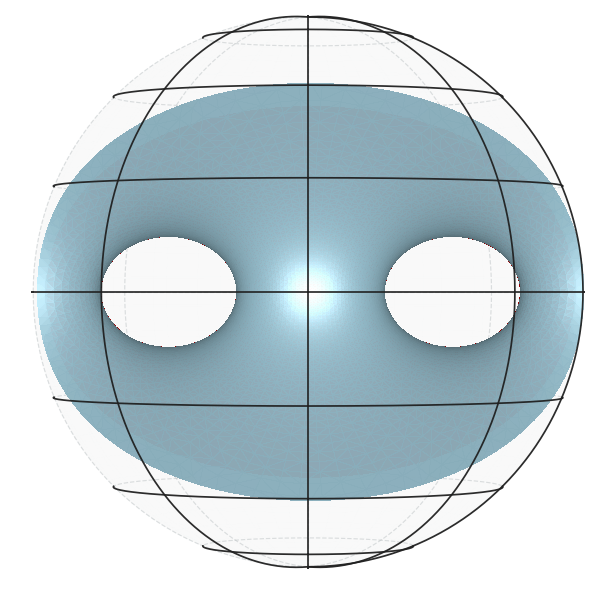}\\
    \includegraphics[width=0.202\textwidth,trim={30 30 14 14},clip]{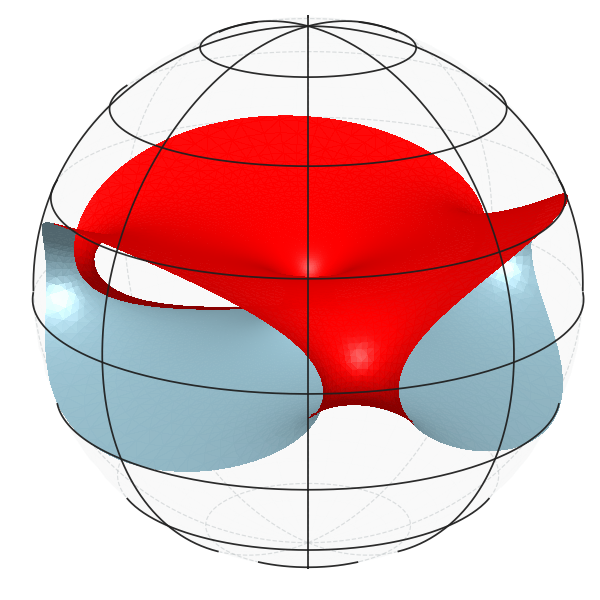}\hspace{2cm}
    \includegraphics[width=0.202\textwidth,trim={30 30 14 14},clip]{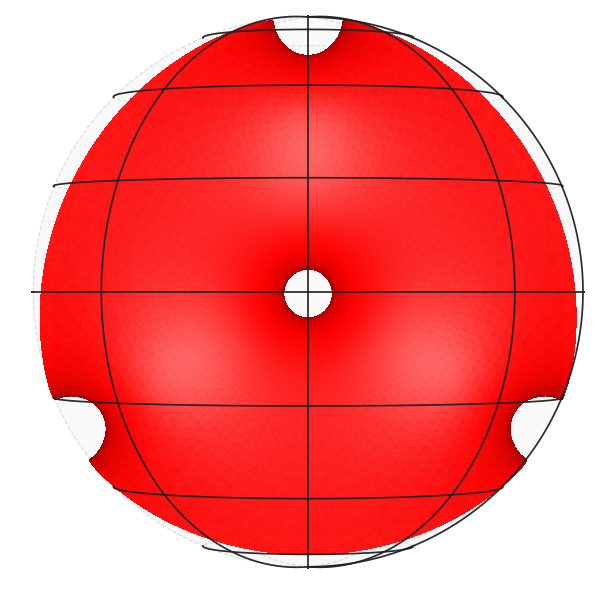}\hspace{2cm}
    \includegraphics[width=0.202\textwidth,trim={30 30 14 14},clip]{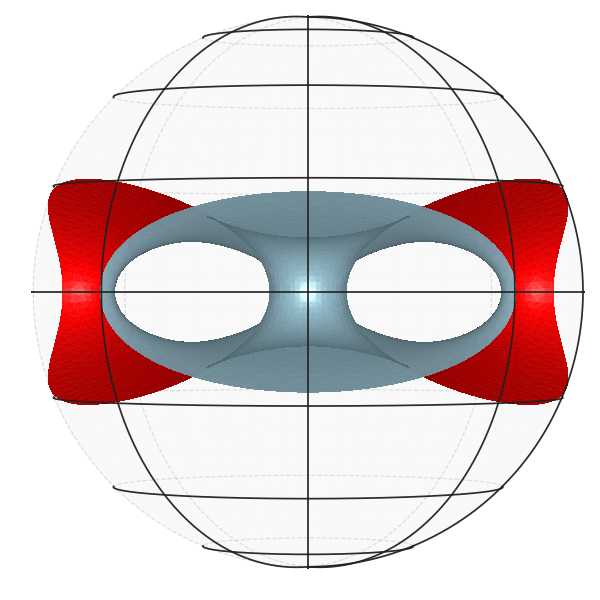}\\
    \includegraphics[width=0.202\textwidth,trim={30 30 14 14},clip]{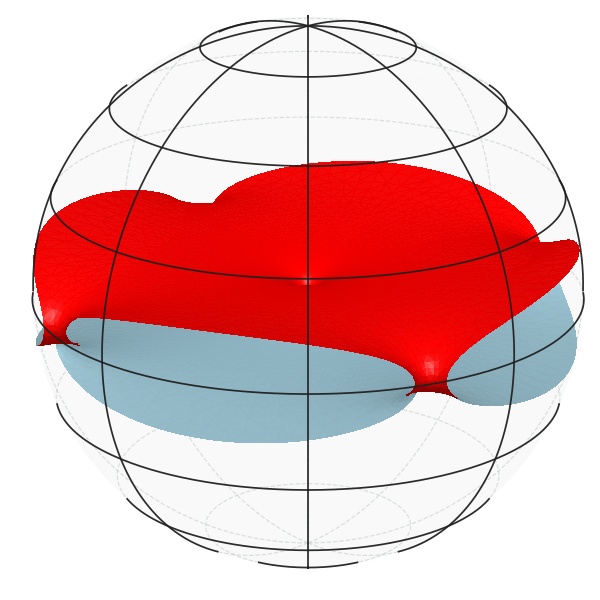}\hspace{2cm}
    \includegraphics[width=0.202\textwidth,trim={30 30 14 14},clip]{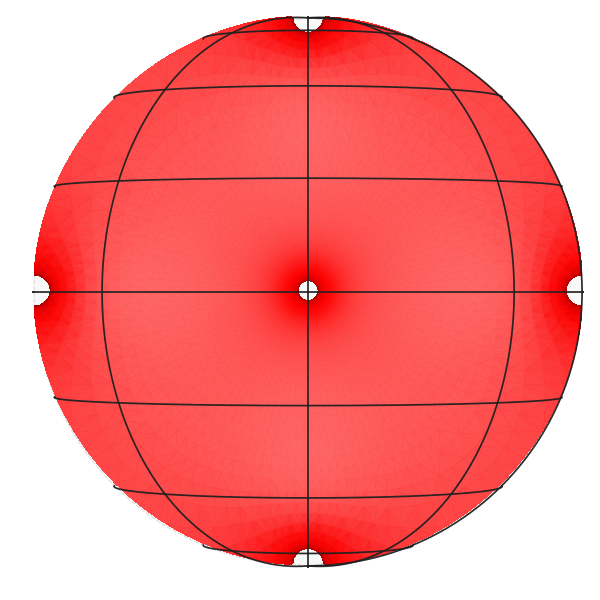}\hspace{2cm}
    \includegraphics[width=0.202\textwidth,trim={30 30 14 14},clip]{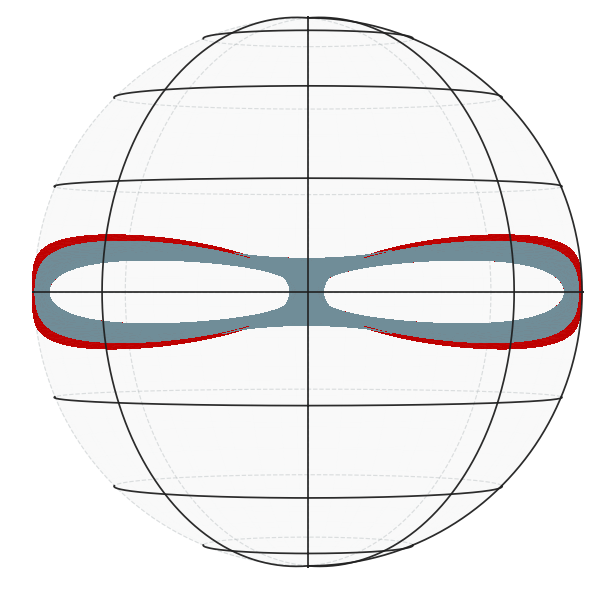}\\
    \includegraphics[width=0.202\textwidth,trim={30 30 14 14},clip]{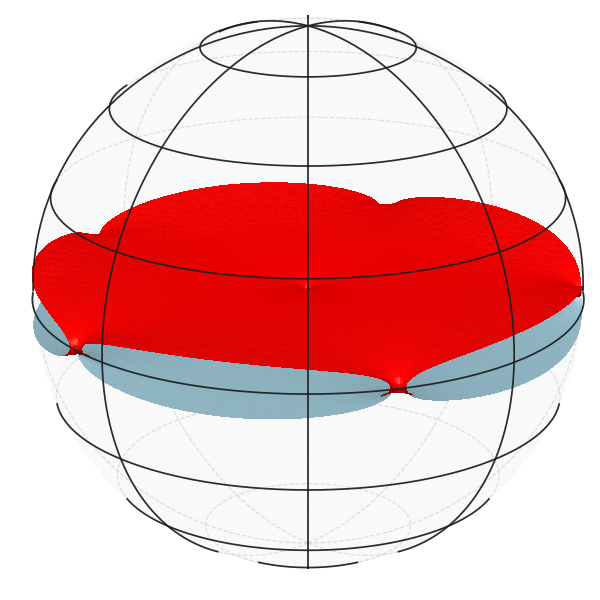}\hspace{2cm}
    \includegraphics[width=0.202\textwidth,trim={30 30 14 14},clip]{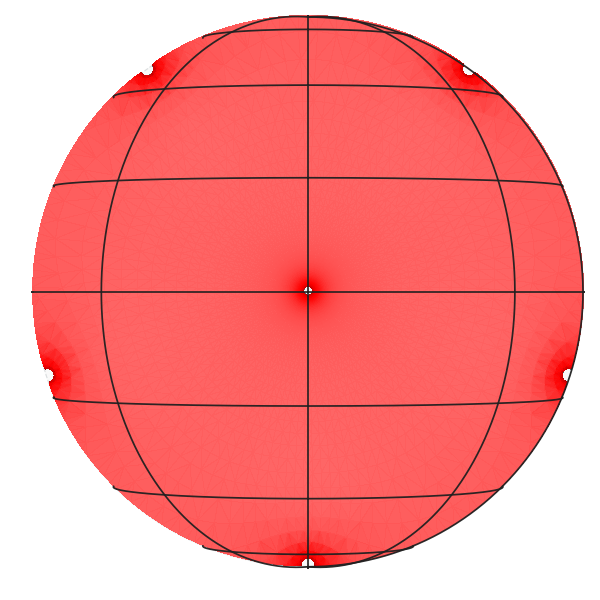}\hspace{2cm}
    \includegraphics[width=0.202\textwidth,trim={30 30 14 14},clip]{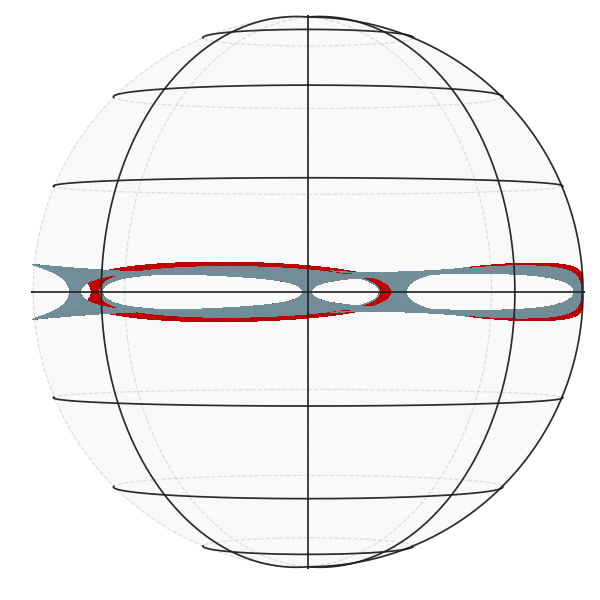}\\
    \includegraphics[width=0.202\textwidth,trim={30 30 14 14},clip]{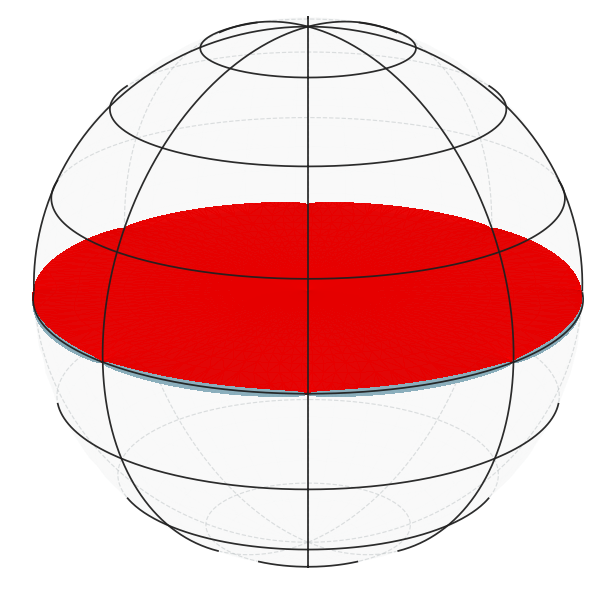}\hspace{2cm}
    \includegraphics[width=0.202\textwidth,trim={30 30 14 14},clip]{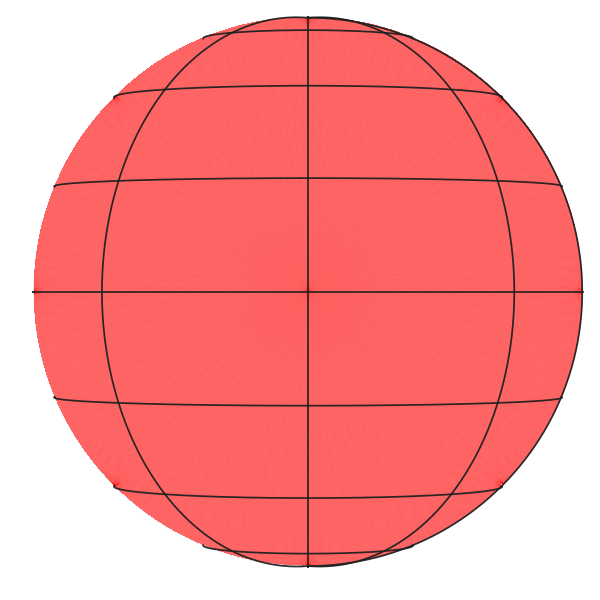}\hspace{2cm}
    \includegraphics[width=0.202\textwidth,trim={30 30 14 14},clip]{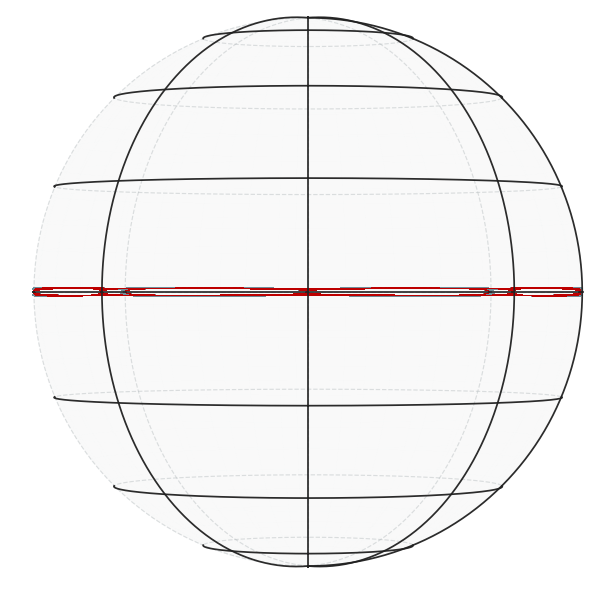}
    \caption{The free boundary minimal tori in the unit ball computed in this paper, with \(b\in\{1,2,3,4,5,8\}\) boundary components (one per row; from left to right an oblique view, the view along the symmetry axis, and a side view), obtained as extremal configurations of the weighted Steklov problem on the torus with \(b\) round holes (Section~\ref{ssec:surfaces}); the moduli \(\tau^\star\) and extremal values are listed in Table~\ref{tab:sanity-checks}.}
  \label{fig:surfaces}
\end{figure}

Our contributions are methodological and computational. On the side of the method: (1)~\emph{Formulation.} We reduce the extremal problem to \eqref{eq:intro-finite-dim-problem}, whose geometric parameters \((\tau,\mathbf c,\mathbf r)\), with one center fixed, cover the \(3b\)-dimensional moduli space of conformal structures (Section~\ref{ssec:circle-domain}). (2)~\emph{Eigensolver.} The objective is evaluated with the doubly periodic layer potentials of \cite{kim2025layer}, for which the weighted Steklov problem is a generalized eigenproblem for a pencil of boundary integral operators, discretized by a spectrally convergent Nystr\"om method (Sections~\ref{ssec:bie} and~\ref{ssec:nystrom}). (3)~\emph{Derivatives.} The left eigenvector of this pencil is the boundary trace of the eigenfunction. Consequently the derivative of \(\bar\sigma_k\) with respect to every component of \(\mathbf p\)---weight, centers, radii, and modulus---is given by a single formula (Proposition~\ref{prop:unified-derivative}) whose discretization is the exact derivative of the discrete eigenvalue; the Hadamard boundary integrals of \cite{Kao_2020} are recovered as special cases, and the modulus derivative, which is new, rests on differentiating the torus Green function through the heat equation for the theta function (Appendix~\ref{sec:analytic formula for dtau}). (4)~\emph{Nonsmooth optimization.} We use gradient ascent while \(\sigma_1\) is simple and, once the first \(m\) eigenvalues coalesce, refine the cluster by maximin steps in the manner of \textcite{Overton_1988} followed by quasi-Newton SQP steps along the active manifold in the manner of \textcite{Bareilles_2023} (Section~\ref{ssec:two-stage-optimization}). A computed critical configuration comes with a certificate of Clarke stationarity, a positive semidefinite matrix \(\Lambda\), and we show that this certificate is the surface (Section~\ref{ssec:stationarity}, Theorem~\ref{prop:certificate}): \(L_\rho\Lambda\) is the rescaling under which the spectral map sends \(\partial\Sigma\) to the unit sphere, and stationarity in the geometric parameters is the conformality of that map, because the parametrization covers the whole moduli space. The optimality conditions of the finite-dimensional problem are thus exactly the free boundary minimal surface conditions. (5)~\emph{Verification.} The immersion obtained from the certificate is verified a posteriori through residuals of the immersion, conformality, and sphere conditions (Section~\ref{ssec:verification}).

On the side of results, we apply the method for \(b\in\{1,2,3,4,5,8\}\) in the class of configurations with the holes on a common closed geodesic and a shared weight profile, and certify stationarity with respect to all parameters, the symmetry-breaking ones included (Table~\ref{tab:sanity-checks}). For \(b=1\) the method finds a critical configuration with \(\bar\sigma_1^\star\approx11.33006466\) and \(\tau^\star\approx i\) at which the first three normalized eigenvalues coalesce; the surface in the first row of Figure~\ref{fig:surfaces} is an embedded genus-one free boundary minimal surface with connected boundary, its coordinate functions (top row of Figure~\ref{fig:eigenfunction_contours}) defining a map into \(\mathbb B^3\) that satisfies the free boundary minimal surface conditions to the accuracy of Table~\ref{tab:sanity-checks}. Whether such a surface exists was the first open question in the survey \cite{li2020}; it was settled by the equivariant min--max construction of \textcite{Carlotto_2022}, which gives existence and a \(\mathbb D_2\) symmetry, and the computation supplies its conformal type, its extremal value, a symmetry group of order eight containing \(\mathbb D_2\), and the fact that its coordinate functions are first Steklov eigenfunctions. For \(b\geq2\), we compute critical configurations with holes equispaced on the geodesic with equal radii and identical weights (Figure~\ref{fig:boundary_densities}), whose surfaces (remaining rows of Figures~\ref{fig:surfaces} and~\ref{fig:eigenfunction_contours}) consist of two nearly flat disks joined by a central catenoidal neck and by \(b\) equatorial necks: the family constructed for large \(b\) by \textcite{Folha_Pacard_Zolotareva_2017} and for \(b=2\) by \textcite{Franz_Ketover_Schulz_2024}, here realized for \(b\in\{2,3,4,5,8\}\) together with their moduli and extremal values; the sheets approach each other as \(b\) grows and are nearly coincident for \(b=8\) (last row of Figure~\ref{fig:surfaces}). Along this family \(\bar\sigma_1^\star\) increases from \(11.3301\) for \(b=1\) to \(12.5661\) for \(b=8\), within \(3\times10^{-4}\) of \(4\pi\approx12.5664\), the value of the doubly covered disk to which the surfaces converge, the modulus \(\Im\tau^\star\) grows, and the smallest eigenvalue of the certificate \(\Lambda^\star\) decreases toward zero as the two sheets approach each other (Section~\ref{ssec:extremal-configurations}).

None of these configurations is a maximizer of \(\bar\sigma_1\), and for \(b\leq5\) not even locally (Section~\ref{ssec:equilateral}). For \(b=1\) the supremum is the Fraser--Schoen bound \(2\pi(\gamma+b)=4\pi\), attained on a family of conformal classes at which the certificate has rank two and the spectral map is a branched double cover of the equatorial disk rather than a surface. \textcite{Carlotto_2022} asked whether their surface maximizes the normalized first Steklov eigenvalue: it is a saddle point, a maximizer within its symmetry class only (Figure~\ref{fig:b1-fixed-camera}). For \(b\geq2\) the largest values we find are at configurations where five or six eigenvalues coalesce, so that the Fraser--Schoen correspondence produces free boundary minimal tori in \(\mathbb B^5\) and \(\mathbb B^6\) rather than \(\mathbb B^3\); that the collinear surfaces are not the global maximizers also follows, for \(b\leq5\), from the genus-zero values of \textcite{Kao_2020} and the gluing inequality of \textcite{Matthiesen_Petrides_2020}. The maximizers are not the geometric objects of interest here, and the surfaces in \(\mathbb B^3\) remain our focus.

\smallskip\noindent\textbf{Outline.} Section~\ref{sec:problem-setup} develops the method at the continuous level, from the conformal reduction and the eigensolver through the derivatives and the optimization to the stationarity certificate and the verification criteria. Section~\ref{sec:implementation} collects discretization and implementation details, Section~\ref{sec:experiments} reports the computations and the resulting surfaces, and we conclude in Section~\ref{sec:discussion} with a discussion. Appendix~\ref{sec:analytic formula for dtau} derives the modulus derivative of the torus Green function and Appendix~\ref{sec:qbx} describes the close-evaluation procedure for the layer potentials.

\subsection*{Notation}

For \(n\in\mathbb N\) we write \([n]:=\{1,\dots,n\}\). Throughout, \(k\) indexes eigenvalues, \(\ell\in[b]\) boundary components, \(n\) Fourier modes, \(j\) components of a vector (of \(\mathbf p\) or of the spectral map), and \(i\) quadrature nodes; \(p\) is a generic component of \(\mathbf p\) and \(\partial_p\) the partial derivative with respect to it; \((\sigma_k,u_k)\) is a Steklov eigenpair and \(\bar\sigma_k\) the length-normalized eigenvalue; \(m\) is the multiplicity of the first eigenvalue at a target configuration; \(N_\ell\) is the number of quadrature nodes on \(\Gamma_\ell\) and \(N:=\sum_{\ell\in[b]}N_\ell\); and a superscript \({}^\star\) marks quantities at a computed critical configuration, e.g.,\ \(\tau^\star\), \(\Omega^\star\), \(\bar\sigma_1^\star\).

\section{Formulation and computational method}
\label{sec:problem-setup}

This section develops the components of the computational method at the continuous level; discretization and implementation details are collected in Section~\ref{sec:implementation}.

\subsection{Extremal Steklov problem and conformal reduction}
\label{ssec:circle-domain}

We make precise the reduction of \eqref{eq:intro-extremal-steklov-problem} to \eqref{eq:intro-finite-dim-problem} sketched in the introduction. Let \(\Sigma\) be a surface of genus \(\gamma=1\) with \(b\) boundary components. By the uniformization theorem for compact Riemann surfaces \cite[Theorem~4.4.1]{Jost2006}, every closed Riemann surface \(S\) of genus one carries a flat metric, unique up to scale. The Koebe uniformization theorem for surfaces of finite type \cite[Theorem~0.2]{HeSchramm1993} then states that the interior of \(\Sigma\) is conformally equivalent to a \emph{circle domain} in such an \(S\), that is, to the complement in \(S\) of finitely many closed disks that are round with respect to the flat metric, and that the pair consisting of \(S\) and the circle domain exists and is unique up to conformal homeomorphism; see \cite{Zeng_2009,Crane_2013} for numerical realizations. Normalizing one lattice period to \(1\), we write \(\mathbb T_\tau:=\mathbb C/\mathcal L_\tau\) with \(\mathcal L_\tau:=\mathbb Z+\tau\mathbb Z\), \(\Im\tau>0\), and call \(\tau\) the conformal modulus. The corresponding domain is the finitely connected flat torus \(\Omega_{\tau,\mathbf c,\mathbf r}:=\mathbb T_\tau\setminus\bigcup_{\ell\in[b]}\overline{\mathbb D_\tau(c_\ell,r_\ell)}\), \(\mathbb D_\tau(c_\ell,r_\ell):=\{z\in\mathbb T_\tau:|z-c_\ell|<r_\ell\}\), the distance being that of the quotient and the removed disks having pairwise disjoint closures. We abbreviate \(\Omega:=\Omega_{\tau,\mathbf c,\mathbf r}\), \(\mathbb D_\ell:=\mathbb D_\tau(c_\ell,r_\ell)\), \(\Gamma_\ell:=\partial\mathbb D_\ell\) when no confusion arises. Thus \(\partial\Omega=\bigcup_{\ell\in[b]}\Gamma_\ell\) is a union of \(b\) circles.

We equip \(\mathbb T_\tau\) and \(\Omega\) with the flat metric \(g_\tau=|dz|^2\). Each boundary component is parametrized by \(\xi_\ell(\theta)=c_\ell+r_\ell e^{i\theta}\), \(0\leq\theta<2\pi\), so that \(ds_\tau=r_\ell\,d\theta\) and the outward conormal of \(\Omega\) is \(\nu_\tau(\xi_\ell(\theta))=-e^{i\theta}\), pointing into the removed disk. This convention fixes the signs below.

The Steklov spectrum of \((\Sigma,g)\) depends on the conformal factor only through its boundary trace, and this is what makes the reduction effective. To see this, let \(\Xi\colon\overline\Omega\to\Sigma\) be the inverse of the uniformizing map, extended smoothly to the boundary, and write \(\Xi^*g=e^{2\omega}g_\tau\) for a smooth real-valued function \(\omega\) on \(\overline\Omega\). Identifying a function on \(\Sigma\) with its pullback \(u\circ\Xi\), we have \(\Delta_{\Xi^*g}u=e^{-2\omega}\Delta_\tau u\), \(\partial_{\nu_{\Xi^*g}}u=e^{-\omega}\partial_{\nu_\tau}u\) and \(ds_{\Xi^*g}=e^{\omega}ds_\tau\), so the geometric Steklov problem \eqref{eq:geometric_steklov_problem} is equivalent to the weighted flat problem
\begin{equation}\label{eq:weighted-steklov-problem}
\Delta_\tau u=0 \quad\text{in }\Omega, \qquad
\partial_{\nu_\tau}u=\sigma\rho u
\quad\text{on }\partial\Omega, \qquad \rho:=e^{\omega}\big|_{\partial\Omega}.
\end{equation}
In particular, \(\sigma_k(\Sigma,g)=\sigma_k(\Omega,\rho)\) and \(L_g(\partial\Sigma)=L_\rho:=\int_{\partial\Omega}\rho\,ds_\tau\), hence, by \eqref{eq:intro-normalized-eigenvalue}, \(\bar\sigma_k(\Sigma,g)=\bar\sigma_k(\Omega,\rho)\). Conversely, every smooth function on \(\partial\Omega\) extends smoothly to \(\overline\Omega\), so optimizing over the boundary trace of \(\omega\) captures every conformal factor relevant to the Steklov spectrum; in the language of \cite[\S2.3.2]{Karpukhin_Metras_2022}, on surfaces the Steklov problem with a boundary density is the classical Steklov problem for a conformally changed metric, and conversely. Since \eqref{eq:weighted-steklov-problem} is real and self-adjoint with respect to the weighted boundary inner product, its eigenfunctions can be taken real and orthonormal in \(L^2(\partial\Omega,\rho\,ds_\tau)\); we fix this normalization throughout:
\begin{equation}\label{eq:weighted-eigenfunction-normalization}
\int_{\partial\Omega}\rho\,u_ku_{k'}\,ds_\tau=\delta_{kk'}.
\end{equation}

For computation, we parametrize the boundary trace of the conformal factor on each circle by a truncated Fourier series: writing \(\omega_\ell:=\omega|_{\Gamma_\ell}\) as a function of the angle \(\theta\),
\begin{equation}
    \omega_\ell(\theta)
    =
    \alpha_{\ell,0}
    +
    \sum_{n=1}^{N_\rho}
    \left[
        \alpha_{\ell,n}\cos(n\theta)
        +
        \beta_{\ell,n}\sin(n\theta)
    \right],
    \qquad
    \rho_\ell:=\rho|_{\Gamma_\ell}=e^{\omega_\ell},
    \qquad \ell\in[b].
\label{eq:fourier-boundary-weight}
\end{equation}
Positivity of the weight is thus automatic; in terms of \eqref{eq:fourier-boundary-weight}, \(L_\rho=\sum_{\ell\in[b]}r_\ell\int_0^{2\pi}e^{\omega_\ell(\theta)}\,d\theta\).

For fixed \(b\), the optimization variables are the modulus \(\tau\), the centers and radii \((\mathbf c,\mathbf r)\) of the disks, and the Fourier coefficients \((\boldsymbol\alpha,\boldsymbol\beta)\) of \(\omega\); we collect them in the parameter vector \(\mathbf p:=(\tau,\mathbf c,\mathbf r,\boldsymbol\alpha,\boldsymbol\beta)\in\mathbb R^{P}\), \(P=2+3b+b(2N_\rho+1)\) as in \eqref{eq:intro-finite-dim-problem}, and write \(\Omega_{\mathbf p}\), \(\rho_{\mathbf p}\), and \(\bar\sigma_k(\mathbf p):=\bar\sigma_k(\Omega_{\mathbf p},\rho_{\mathbf p})\) when the dependence on the parameters should be emphasized. The extremal problem \eqref{eq:intro-extremal-steklov-problem} for genus one becomes the parameterized problem \eqref{eq:intro-finite-dim-problem}, with the constraints made explicit:
\begin{equation}\label{eq:objective}
\begin{aligned}
\max_{\mathbf p}\ \bar\sigma_1(\mathbf p)\quad\text{subject to}\quad &\Im\tau>0,\quad r_\ell>0,\quad \operatorname{dist}_{\mathbb T_\tau}(c_\ell,c_{\ell'})>r_\ell+r_{\ell'}\ \ \text{for }1\leq\ell<\ell'\leq b,\\
&c_1=\bar c_1,\quad \alpha_{1,0}=0,
\end{aligned}
\end{equation}
where \(\operatorname{dist}_{\mathbb T_\tau}\) denotes the flat distance on the torus and \(\bar c_1\) is a prescribed point of \(\mathbb T_\tau\); the last two constraints are normalizations, explained below. The weighted Steklov problem \eqref{eq:weighted-steklov-problem} is an implicit constraint in \eqref{eq:objective}: \(\sigma_1(\mathbf p)\) is the first positive eigenvalue of \eqref{eq:weighted-steklov-problem} on \(\Omega_{\mathbf p}\) with weight \(\rho_{\mathbf p}\).

The parametrization is complete. With \(c_1\) fixed, the geometric parameters \((\tau,\mathbf c,\mathbf r)\) have \(3b\) real components, the dimension \(6\gamma-6+3b\) of the moduli space of genus-one surfaces with \(b\) boundary components, and by the existence and uniqueness of the circle-domain representation \cite[Theorem~0.2]{HeSchramm1993} they cover it; \(\omega\) supplies the conformal factor, and \(\alpha_{1,0}=0\) normalizes it, since the shift \(\omega\mapsto\omega+\text{const}\) rescales \(\rho\) and \(L_\rho\) by a common factor and \(\sigma_k\) by its inverse, leaving \(\bar\sigma_1\) unchanged. The remaining redundancies, the modular action on \(\tau\) (\(\tau\sim\tau+1\sim-1/\tau\)) and the relabeling of the holes, are discrete and irrelevant for a local method; the reported configurations are representatives of their orbits.

\subsection{Boundary-integral formulation}
\label{ssec:bie}

We use the doubly periodic layer-potential formulation of \textcite{kim2025layer} as the forward eigensolver and summarize it here. Since \(\mathbb T_\tau\) is compact, it admits no Green function in the classical sense: the left-hand side of \(-\Delta_\tau G=\delta_0\) integrates to zero over the torus while the right-hand side does not. One works instead with a modified Green function \(G_\tau\), a solution, unique up to an additive constant, of
\[-\Delta_\tau G_\tau=\delta_0-\frac{1}{|\mathbb T_\tau|},\qquad |\mathbb T_\tau|=\Im\tau,\]
whose right-hand side has zero mean to match the left. It is given explicitly by
\begin{equation}\label{eq:torus-green-function}
G_\tau(z)=-\frac{1}{2\pi} \log|\vartheta_1(z;\tau)|
+\frac{(\Im z)^2}{2\Im\tau},
\end{equation}
where \(\vartheta_1(z;\tau)\) is the Jacobi theta function
\[\vartheta_1(z;\tau)
:=-i\sum_{n\in\mathbb Z}(-1)^n
e^{\pi i\tau(n+1/2)^2}e^{(2n+1)\pi iz}.\]
Any other solution differs from \eqref{eq:torus-green-function} by a constant, and the choice is arbitrary below because \(G_\tau\) is only ever integrated against densities of zero mean.\footnote{The mean of \eqref{eq:torus-green-function} over \(\mathbb T_\tau\) is \(-\tfrac{1}{2\pi}\log|\eta(\tau)|\), with \(\eta\) the Dedekind eta function; the zero-mean solution is obtained by replacing \(\vartheta_1(z;\tau)\) with \(\vartheta_1(z;\tau)/\eta(\tau)\).} For a boundary density \(\varphi\), define the single-layer potential and the adjoint double-layer potential
\begin{subequations}\label{eq:layer-potentials}
\begin{align}
\mathcal S[\varphi](z)
&:=\int_{\partial\Omega}G_\tau(z-\xi)\varphi(\xi)\,ds_\tau(\xi),
&&z\in\overline\Omega,
\label{eq:single-layer-potential}\\
\mathcal D^*[\varphi](z)
&:=\operatorname{p.v.}\!\int_{\partial\Omega}
\partial_{\nu_{\tau}(z)}G_\tau(z-\xi)\varphi(\xi)\,ds_\tau(\xi),
&&z\in\partial\Omega,
\label{eq:adjoint-double-layer-potential}
\end{align}
\end{subequations}
where \(\partial_{\nu_\tau(z)}G_\tau(z-\xi)=
\nabla_zG_\tau(z-\xi)\cdot\nu_\tau(z)\); the single-layer potential is continuous across \(\partial\Omega\), and the principal value in \eqref{eq:adjoint-double-layer-potential} is not needed on a smooth boundary, where the kernel is bounded. The modified single-layer potential is given by
\begin{equation}\label{eq:modified-single-layer-potential}
\mathcal S_{0}[\varphi]
:=(\mathcal S(I-M)+M)[\varphi], \qquad M[\varphi]:=\frac{1}{|\partial\Omega|}\int_{\partial \Omega}\varphi(\xi)\,ds_\tau(\xi).
\end{equation}
Here \(|\partial\Omega|=\int_{\partial\Omega}ds_\tau=2\pi\sum_{\ell\in[b]}r_\ell\) is the flat length of the boundary. It is essential that \eqref{eq:modified-single-layer-potential} uses the unweighted
flat arclength measure, the same measure as in \eqref{eq:single-layer-potential}: it makes the density entering \(\mathcal S\) have zero total mass, so that the background term in
\eqref{eq:torus-green-function} drops out and \(\mathcal S_0[\varphi]\) is harmonic in \(\Omega\).  The Steklov weight \(\rho\) enters the
generalized eigenproblem separately. \textcite{kim2025layer} characterize the eigenpairs of \eqref{eq:weighted-steklov-problem} in terms of these operators.

\begin{theorem}[Steklov boundary-integral formulation \cite{kim2025layer}] \label{thm:weighted-steklov-bie}
Fix a circle domain \(\Omega\subset\mathbb T_\tau\) as in Section~\ref{ssec:circle-domain} and a smooth
positive boundary weight \(\rho\).  A pair \((\sigma_k,u_k)\) is an eigenpair of
\eqref{eq:weighted-steklov-problem} if and only if there is a nonzero
boundary density \(\varphi_k\) satisfying
\begin{equation}\label{eq:weighted-steklov-bie}
\mathsf A[\varphi_k]=\sigma_k\,\mathsf B[\varphi_k]\quad\text{on }\partial\Omega,\qquad
\mathsf A:=(\mathcal D^*+\tfrac12I)(I-M),\quad \mathsf B:=\rho\,\mathcal S_{0}.
\end{equation}
The associated eigenfunction is
\begin{equation}\label{eq:eigenfunction-layer-representation}
u_k(z)=\mathcal S_{0}[\varphi_k](z),
\qquad z\in\Omega.\end{equation}\end{theorem}
The Nystr\"om discretization of \eqref{eq:weighted-steklov-bie}, and the evaluation of \eqref{eq:eigenfunction-layer-representation} and its derivatives in \(\overline\Omega\), are described in Section~\ref{sec:implementation}.

\subsection{Eigenvalue derivatives}
\label{ssec:analytic-sensitivities}

We now derive the derivatives of \(\bar\sigma_k=L_\rho\sigma_k\) with respect to all components of the parameter vector \(\mathbf p\). Throughout, \((\sigma_k,u_k)\) is a \emph{simple} real eigenpair of \eqref{eq:weighted-steklov-problem}, normalized according to \eqref{eq:weighted-eigenfunction-normalization}, and \(\varphi_k\) is its density in Theorem~\ref{thm:weighted-steklov-bie}; clusters of eigenvalues are treated at the end of the section. The point of view is that \eqref{eq:weighted-steklov-bie} is a generalized eigenvalue problem for the pencil \((\mathsf A(\mathbf p),\mathsf B(\mathbf p))\) of boundary operators, which depends smoothly on every parameter, and that the left eigenvector of this pencil is the boundary trace of \(u_k\). All derivatives then follow from one formula, with no adjoint problem and no interior evaluation of \(u_k\) or its gradient.

To make the parameter dependence explicit we represent a function on \(\partial\Omega\) as a function of \((\ell,\theta)\) through the parametrization \(\xi_\ell(\theta)=c_\ell+r_\ell e^{i\theta}\) of Section~\ref{ssec:circle-domain}, so that \(\mathsf A(\mathbf p)\) and \(\mathsf B(\mathbf p)\) act on a fixed space of functions of \((\ell,\theta)\). For a kernel function \(\kappa(w;\theta,\theta')\), which may depend on the target angle \(\theta\) and the source angle \(\theta'\), let
\begin{equation}\label{eq:kernel-operators}
\bigl(\operatorname{op}(\kappa)\varphi\bigr)(\ell,\theta):=\sum_{\ell'\in[b]}\int_0^{2\pi}\kappa\bigl(\xi_{\ell}(\theta)-\xi_{\ell'}(\theta');\theta,\theta'\bigr)\,\varphi(\ell',\theta')\,r_{\ell'}\,d\theta'.
\end{equation}
In this representation \(\mathcal S=\operatorname{op}(\kappa_{\mathcal S})\) and \(\mathcal D^*=\operatorname{op}(\kappa_{\mathcal D^*})\) with \(\kappa_{\mathcal S}(w):=G_\tau(w)\) and \(\kappa_{\mathcal D^*}(w;\theta):=-\nabla G_\tau(w)\cdot e^{i\theta}\), since \(-e^{i\theta}=\nu_\tau(\xi_{\ell}(\theta))\), and \(M=\operatorname{op}(|\partial\Omega|^{-1})\) is the operator with constant kernel. The dependence on \(\mathbf p\) is now visible: \(\tau\) enters through the kernel \(G_\tau\); \(\mathbf c\) and \(\mathbf r\) through the kernel arguments \(\xi_{\ell}(\theta)-\xi_{\ell'}(\theta')\), the weights \(r_{\ell'}\,d\theta'\) and \(|\partial\Omega|\); and \((\boldsymbol\alpha,\boldsymbol\beta)\) through \(\rho\), which appears only in \(\mathsf B\). We write \(\mathbf 1_\ell\) for multiplication by the indicator of \(\Gamma_\ell\), so that \(\mathbf 1_\ell\mathcal Q\) restricts the target and \(\mathcal Q\mathbf 1_\ell\) the source of an operator \(\mathcal Q\) to \(\Gamma_\ell\), and \([\mathbf 1_\ell,\mathcal Q]:=\mathbf 1_\ell\mathcal Q-\mathcal Q\mathbf 1_\ell\). For the complex parameter \(\tau=\tau_1+i\tau_2\) we use the Wirtinger derivative \(\partial_\tau:=\tfrac12(\partial_{\tau_1}-i\partial_{\tau_2})\), so that \(\partial_{\tau_1}\bar\sigma_k=2\Re(\partial_\tau\bar\sigma_k)\), \(\partial_{\tau_2}\bar\sigma_k=-2\Im(\partial_\tau\bar\sigma_k)\). Finally, \(\dot{\mathcal S}:=\operatorname{op}(\partial_\tau \kappa_{\mathcal S})\) and \(\dot{\mathcal D}^*:=\operatorname{op}(\partial_\tau \kappa_{\mathcal D^*})\) denote the operators with differentiated kernels.

\begin{proposition}[Derivatives through the boundary operators]\label{prop:unified-derivative}
For every component \(p\) of \(\mathbf p\),
\begin{equation}\label{eq:unified-derivative}
\partial_p\sigma_k=\int_{\partial\Omega}u_k\,\bigl[(\partial_p\mathsf A-\sigma_k\,\partial_p\mathsf B)\varphi_k\bigr]\,ds_\tau,\qquad
\partial_p\bar\sigma_k=L_\rho\,\partial_p\sigma_k+\sigma_k\,\partial_pL_\rho,
\end{equation}
where \(\partial_p\mathsf A\), \(\partial_p\mathsf B\) and \(\partial_pL_\rho\) are tabulated below (the row for \(c_\ell\) is \(\mathbb R^2\)-valued and gives \(\nabla_{c_\ell}\)):
\begin{center}\footnotesize
\renewcommand{\arraystretch}{1.7}
\begin{tabular}{@{}llll@{}}
\toprule
\(p\) & \(\partial_p\mathsf A\) & \(\partial_p\mathsf B\) & \(\partial_pL_\rho\)\\
\midrule
\(\alpha_{\ell,n}\) & \(0\) & \(\rho_\ell(\theta)\cos(n\theta)\,\mathbf 1_\ell\,\mathcal S_0\) & \(r_\ell\int_0^{2\pi}\rho_\ell(\theta)\cos(n\theta)\,d\theta\)\\
\(\beta_{\ell,n}\) & \(0\) & \(\rho_\ell(\theta)\sin(n\theta)\,\mathbf 1_\ell\,\mathcal S_0\) & \(r_\ell\int_0^{2\pi}\rho_\ell(\theta)\sin(n\theta)\,d\theta\)\\
\(\tau\) & \(\dot{\mathcal D}^*(I-M)\) & \(\rho\,\dot{\mathcal S}(I-M)\) & \(0\)\\
\(c_\ell\in\mathbb R^2\) & \(\bigl[\mathbf 1_\ell,\operatorname{op}(\nabla_w\kappa_{\mathcal D^*})\bigr](I-M)\) & \(\rho\,\bigl[\mathbf 1_\ell,\operatorname{op}(\nabla_w\kappa_{\mathcal S})\bigr](I-M)\) & \(0\)\\
\(r_\ell\) & \((\partial_{r_\ell}\mathcal D^*)(I-M)-(\mathcal D^*+\tfrac12I)\,\partial_{r_\ell}M\) & \(\rho\,\bigl[(\partial_{r_\ell}\mathcal S)(I-M)+(I-\mathcal S)\,\partial_{r_\ell}M\bigr]\) & \(\int_0^{2\pi}\rho_\ell(\theta)\,d\theta\)\\
\bottomrule
\end{tabular}
\end{center}
with, for \(\mathcal Q\in\{\mathcal S,\mathcal D^*\}\) with kernel \(\kappa\),
\begin{subequations}\label{eq:radius-operator-derivatives}
\begin{align}
\partial_{r_\ell}\mathcal Q&=\mathbf 1_\ell\,\operatorname{op}\bigl(e^{i\theta}\!\cdot\nabla_w\kappa\bigr)-\operatorname{op}\bigl(e^{i\theta'}\!\cdot\nabla_w\kappa\bigr)\,\mathbf 1_\ell+\mathcal Q\,\frac{\mathbf 1_\ell}{r_\ell},\label{eq:radius-operator-derivative-K}\\
\partial_{r_\ell}M&=M\Bigl(\frac{\mathbf 1_\ell}{r_\ell}-\frac{2\pi}{|\partial\Omega|}\Bigr).\label{eq:radius-operator-derivative-M}
\end{align}
\end{subequations}
\end{proposition}

\begin{proof}
Let \(\langle f,g\rangle:=\int_{\partial\Omega}fg\,ds_\tau\) be the bilinear pairing at the base point \(\mathbf p\). For any two densities \(\varphi,\chi\), the functions \(\mathcal S_0\varphi\) and \(\mathcal S_0\chi\) are harmonic in \(\Omega\) with normal derivatives \(\mathsf A\varphi\) and \(\mathsf A\chi\), so Green's second identity gives \(\langle\mathcal S_0\chi,\mathsf A\varphi\rangle=\langle\mathsf A\chi,\mathcal S_0\varphi\rangle\), i.e.,\ \(\mathsf A^{\mathsf T}\mathcal S_0=\mathcal S_0^{\mathsf T}\mathsf A\) for the transposes with respect to \(\langle\cdot,\cdot\rangle\) (all kernels are real, so on real functions these are the \(L^2\) adjoints; the bilinear pairing is used because the derivative in \(\tau\) is complex). Consequently the boundary trace \(\psi_k:=\mathcal S_0\varphi_k=u_k|_{\partial\Omega}\) satisfies \(\mathsf A^{\mathsf T}\psi_k=\mathcal S_0^{\mathsf T}\mathsf A\varphi_k=\sigma_k\mathcal S_0^{\mathsf T}\rho\,\mathcal S_0\varphi_k=\sigma_k\mathsf B^{\mathsf T}\psi_k\): the left eigenvector of the pencil is the trace of \(u_k\), and \(\langle\psi_k,\mathsf B\varphi_k\rangle=\int_{\partial\Omega}\rho\,u_k^2\,ds_\tau=1\) by \eqref{eq:weighted-eigenfunction-normalization}. The identity \(\mathsf A(\mathbf p)\varphi_k(\mathbf p)=\sigma_k(\mathbf p)\mathsf B(\mathbf p)\varphi_k(\mathbf p)\) is an identity between functions of \((\ell,\theta)\) and may be differentiated in \(p\) pointwise (all terms are linear in the perturbation, so for \(\tau\) the Wirtinger derivative may be used). Pairing the result with \(\psi_k\), the terms containing \(\partial_p\varphi_k\) cancel, since \(\langle\psi_k,\mathsf A\partial_p\varphi_k\rangle=\langle\mathsf A^{\mathsf T}\psi_k,\partial_p\varphi_k\rangle=\sigma_k\langle\psi_k,\mathsf B\partial_p\varphi_k\rangle\), and the first identity in \eqref{eq:unified-derivative} follows; the second is the product rule, and the expressions for \(\partial_pL_\rho\) are obtained by differentiating \(L_\rho=\sum_{\ell\in[b]}r_\ell\int_0^{2\pi}e^{\omega_\ell(\theta)}\,d\theta\).

The entries \(\partial_p\mathsf A\), \(\partial_p\mathsf B\) are obtained by differentiating \eqref{eq:kernel-operators} and \(M\) under the integral signs. For \(\alpha_{\ell,n}\), \(\beta_{\ell,n}\) only \(\rho=e^{\omega}\) depends on the parameter, with \(\partial_{\alpha_{\ell,n}}\rho=\rho_\ell\cos(n\theta)\mathbf 1_\ell\). For \(\tau\) only the kernel depends on the parameter, and \(M\) does not. Since \(G_\tau\) is real but not holomorphic in \(\tau\), \(\partial_\tau G_\tau\) is complex-valued; the logarithmic singularity of \(G_\tau\) does not depend on \(\tau\), so \(\partial_\tau G_\tau\) is smooth and no principal value is needed in \(\dot{\mathcal D}^*\). Explicitly,
\begin{equation}\label{eq:dtau-G-closed}
\partial_\tau G_\tau(z)=-\frac{1}{4\pi}\,\frac{\partial_\tau\vartheta_1(z;\tau)}{\vartheta_1(z;\tau)}+\frac{i(\Im z)^2}{4(\Im\tau)^2},
\end{equation}
which is derived in Appendix~\ref{sec:analytic formula for dtau} together with the series \eqref{eq:K_series_impl} by which it is evaluated and its expression \eqref{eq:AppA_final_dtau_G} through Weierstrass functions. A translation of \(\Gamma_\ell\) by a vector \(v\) changes the kernel argument \(w=\xi_{\ell}(\theta)-\xi_{\ell'}(\theta')\) by \(+v\) when the target lies on \(\Gamma_\ell\) and by \(-v\) when the source does, and leaves \(\nu_\tau\), the weights, and \(M\) unchanged; the derivative along \(v\) is therefore \([\mathbf 1_\ell,\operatorname{op}(\nabla_w\kappa\cdot v)]\), and the gradient is the commutator with the vector-valued kernel \(\nabla_w\kappa\). A change of \(r_\ell\) moves the target point by \(e^{i\theta}\,dr_\ell\) and the source point by \(e^{i\theta'}\,dr_\ell\) when they lie on \(\Gamma_\ell\), leaves \(\nu_\tau\) unchanged, scales the weight \(r_\ell\,d\theta'\) of sources on \(\Gamma_\ell\), and changes \(|\partial\Omega|\) in \(M\); this gives \eqref{eq:radius-operator-derivative-K}, whose first two terms have the structure of the commutator with \(v\) replaced by the radial directions on the target and source sides; \eqref{eq:radius-operator-derivative-M} is \eqref{eq:radius-operator-derivative-K} applied to the constant kernel of \(M\) together with \(\partial_{r_\ell}|\partial\Omega|=2\pi\), and the entries for \(\mathsf A=(\mathcal D^*+\tfrac12I)(I-M)\) and \(\mathsf B=\rho(\mathcal S(I-M)+M)\) follow from the product rule.
\end{proof}

We comment on the individual rows of the table in Proposition~\ref{prop:unified-derivative} and on their relation to the Hadamard-type formulas of \cite{Kao_2020}.

\emph{Weight.} Since \(\mathcal S_0\varphi_k=u_k\), the first row gives \(\partial_{\alpha_{\ell,n}}\sigma_k=-\sigma_k\int_{\Gamma_\ell}\rho_\ell\cos(n\theta)\,u_k^2\,ds_\tau\). More generally, for a variation \(\delta\rho\) of the weight,
\begin{equation}\label{eq:weight-sensitivities}
\delta_\rho\sigma_k=-\sigma_k\int_{\partial\Omega}u_k^2\,\delta\rho\,ds_\tau,\qquad
\delta_\rho\bar\sigma_k=\sigma_k\int_{\partial\Omega}\bigl(1-L_\rho u_k^2\bigr)\,\delta\rho\,ds_\tau,
\end{equation}
which is \cite[Proposition~2.6]{Kao_2020}; with \(\rho=e^\omega\) a variation \(\delta\omega\) of the log-weight induces \(\delta\rho=\rho\,\delta\omega\), and for the parametrization \eqref{eq:fourier-boundary-weight},
\[\partial_{\alpha_{\ell,n}}\bar\sigma_k=\sigma_k\int_{\Gamma_\ell}\bigl(1-L_\rho u_k^2\bigr)\rho_\ell\cos(n\theta)\,ds_\tau,\qquad
\partial_{\beta_{\ell,n}}\bar\sigma_k=\sigma_k\int_{\Gamma_\ell}\bigl(1-L_\rho u_k^2\bigr)\rho_\ell\sin(n\theta)\,ds_\tau .\]

\emph{Modulus.} The third row gives \(\partial_\tau\bar\sigma_k=L_\rho\int_{\partial\Omega}u_k\,[\dot{\mathcal D}^*(I-M)-\sigma_k\rho\,\dot{\mathcal S}(I-M)]\varphi_k\,ds_\tau\), with the kernel \eqref{eq:dtau-G-closed}. The derivative is taken with the centers \(c_\ell\in\mathbb C\) fixed in the lift in which they are given; a different representative of \(c_\ell\) modulo \(\mathcal L_\tau\) changes \(\partial_\tau\bar\sigma_k\) by a combination of the center derivatives, and the same choice is implicit in \(\partial_\tau G_\tau(z-\xi)\); at a stationary point, where the center derivatives vanish, the modulus derivative is independent of the lift. This entry has no counterpart in \cite{Kao_2020}, where the planar domain has no modulus; the derivative with respect to the conformal class is the ingredient identified in \cite{Kao_2022} as a challenge for computations beyond genus zero.

\emph{Centers.} In the commutator \([\mathbf 1_\ell,\operatorname{op}(\nabla_w\kappa)]\) the self-interaction block of \(\Gamma_\ell\) cancels, so only hole-to-hole blocks enter and every kernel involved, \(\nabla_wG_\tau\) for \(\mathcal S\) and the Hessian of \(G_\tau\) for \(\mathcal D^*\), is evaluated at positive distance and is smooth. For \(b=1\) the derivative vanishes identically, as it must by translation invariance.

\emph{Radii.} The self-interaction block of \(\partial_{r_\ell}\mathcal Q\) involves the kernel at \(w=r_\ell(e^{i\theta}-e^{i\theta'})\); the logarithmic part of \(G_\tau\) contributes to it only explicit smooth terms (Section~\ref{ssec:nystrom}), and the term \(\mathcal Q\,\mathbf 1_\ell/r_\ell\) needs nothing new.

\emph{Relation to the Hadamard formulas.} For the radii and centers, the table agrees with the Hadamard shape derivative of a simple weighted Steklov eigenvalue \cite[Proposition~2.10]{Kao_2020}, taken with the weight transported along the moving boundary,
\begin{subequations}\label{eq:hadamard-derivatives}
\begin{align}
\partial_{r_\ell}\bar\sigma_k&=\int_{\Gamma_\ell}\Bigl[\bigl(1-L_\rho u_k^2\bigr)\frac{\sigma_k\rho_\ell}{r_\ell}-L_\rho|\nabla u_k|^2+2L_\rho\rho_\ell^2\sigma_k^2u_k^2\Bigr]ds_\tau,\label{eq:dtilde_drj}\\
\nabla_{c_\ell}\bar\sigma_k&=-L_\rho\int_{\Gamma_\ell}\Bigl[F_k\,e_{r,\ell}(\theta)-\frac{\sigma_k}{r_\ell}(\partial_\theta\rho_\ell)\,u_k^2\,e_{\theta,\ell}(\theta)\Bigr]ds_\tau,\label{eq:dtilde_dcj}
\end{align}
\end{subequations}
with \(e_{r,\ell}(\theta):=(\cos\theta,\sin\theta)\), \(e_{\theta,\ell}(\theta):=(-\sin\theta,\cos\theta)\), \(F_k:=|\nabla u_k|^2-2\rho_\ell^2\sigma_k^2u_k^2+\sigma_k\rho_\ell u_k^2/r_\ell\), \(\partial_\theta\rho_\ell=\rho_\ell\,\partial_\theta\omega_\ell\), and \(|\nabla u_k|^2=(r_\ell^{-1}\partial_\theta u_k)^2+(\sigma_k\rho_\ell u_k)^2\) on \(\Gamma_\ell\) by the Steklov condition, so that \eqref{eq:dtilde_drj} and \eqref{eq:dtilde_dcj} involve only the trace of \(u_k\) and its tangential derivative.
For the computations we use Proposition~\ref{prop:unified-derivative}: it covers all parameters with one formula; it needs only the trace of \(u_k\), with no curvature terms or orientation conventions; and applied to the Nyström matrices \(\mathsf A_h(\mathbf p)\), \(\mathsf B_h(\mathbf p)\) it gives the exact derivative of the discrete eigenvalue, so the gradient the optimizer uses is consistent with the objective.

\smallskip\noindent\textbf{Derivatives of eigenvalue clusters.} Stage~2 of Section~\ref{ssec:two-stage-optimization} acts on a cluster \(\{\sigma_k\}_{k\in[m]}\) of \(m\) nearly coalesced eigenvalues and needs, in place of the individual derivatives, the derivative of the cluster as a whole. As long as the cluster stays separated from the rest of the spectrum, the pencil \((\mathsf A,\mathsf B)\) has smooth bases \(\Phi(\mathbf p)\), \(\Psi(\mathbf p):=\mathcal S_0\Phi(\mathbf p)\) of the right and left invariant subspaces belonging to the cluster, normalized by \(\Psi^{\mathsf T}\mathsf B\Phi=I\), and the cluster eigenvalues are the eigenvalues of the smooth symmetric \(m\times m\) matrix functions
\begin{equation}\label{eq:cluster-matrix}
\mathsf S(\mathbf p):=\Psi^{\mathsf T}\mathsf A\Phi,\qquad \bar{\mathsf S}(\mathbf p):=L_\rho\,\mathsf S(\mathbf p).
\end{equation}
At an exactly \(m\)-fold eigenvalue \(\sigma\), differentiating \eqref{eq:cluster-matrix} and using \(\mathsf A\Phi=\sigma\mathsf B\Phi\), \(\Psi^{\mathsf T}\mathsf A=\sigma\Psi^{\mathsf T}\mathsf B\) and the derivative of the normalization \(\Psi^{\mathsf T}\mathsf B\Phi=I\) shows that the terms containing \(\partial_p\Phi\) and \(\partial_p\Psi\) cancel, exactly as in the proof of Proposition~\ref{prop:unified-derivative}, and
\begin{equation}\label{eq:cluster-derivative}
\partial_p\mathsf S=\Psi^{\mathsf T}\bigl(\partial_p\mathsf A-\sigma\,\partial_p\mathsf B\bigr)\Phi,\qquad
\partial_p\bar{\mathsf S}=L_\rho\,\partial_p\mathsf S+(\partial_pL_\rho)\,\sigma I.
\end{equation}
By the proof of Proposition~\ref{prop:unified-derivative}, \(\Psi\) consists of the boundary traces of the \(u_k\), and the entries of \eqref{eq:cluster-derivative} are the polarized forms of the individual derivatives, e.g.,\ \(\bigl(\delta_\rho\bar{\mathsf S}\bigr)_{kk'}=\sigma\int_{\partial\Omega}(\delta_{kk'}-L_\rho u_ku_{k'})\,\delta\rho\,ds_\tau\) for the weight, whose diagonal is the second formula in \eqref{eq:weight-sensitivities}. The implementation uses \eqref{eq:cluster-derivative} with the discrete pencil in place of the continuous one.

\subsection{Two-stage optimization procedure}
\label{ssec:two-stage-optimization}

The objective \(\bar\sigma_1(\mathbf p)\) is smooth wherever \(\sigma_1(\mathbf p)\) is simple, with gradient given by Proposition~\ref{prop:unified-derivative}, but it is only Lipschitz where \(\sigma_1\) is multiple, and multiplicity is exactly what we expect at a maximizer: by the Fraser--Schoen correspondence, an extremal metric whose first eigenspace has dimension \(m\) produces a spectral map into \(\mathbb B^m\). The target configurations therefore satisfy
\begin{equation}\label{eq:intro-three-dimensional-cluster}
0=\bar\sigma_0<\bar\sigma_1=\dots=\bar\sigma_m<\bar\sigma_{m+1}
\end{equation}
for some \(m\geq2\) (a harmonic map into \([-1,1]\) with \(|X|=1\) on the boundary is constant, so \(m=1\) is impossible); the surfaces in \(\mathbb B^3\) of Sections~\ref{ssec:extremal-configurations}--\ref{ssec:surfaces} have \(m=3\), and the maximizers of Section~\ref{ssec:equilateral} have clusters of size up to six. This is the classical setting of nonsmooth eigenvalue optimization \cite{Overton_1988,Lewis_Overton_1996}. In terms of the cluster matrix \eqref{eq:cluster-matrix}, \(\bar\sigma_1=\lambda_{\min}\circ\bar{\mathsf S}\) is the composition of a smooth map with the concave, nonsmooth function \(\lambda_{\min}\), and at a point of exact multiplicity the derivative of \(\bar{\mathsf S}\) is \eqref{eq:cluster-derivative}. The set of parameters at which the cluster is an exact \(m\)-fold eigenvalue is, generically, a smooth manifold of codimension \(\tfrac12m(m+1)-1\), on which \(\bar\sigma_1\) is smooth; this is the active manifold in the sense of \cite{Bareilles_2023}.

\emph{Stage 1: gradient ascent.} While the first eigenvalue is simple and numerically separated from \(\bar\sigma_2\), we perform gradient ascent on \(\bar\sigma_1\) using Proposition~\ref{prop:unified-derivative} with \(k=1\); steps that violate the feasibility conditions of \eqref{eq:objective} are shortened. At every iterate \(\mathbf p\) the eigenproblem is solved anew and the normalized spectrum is ordered as \(0=\bar\sigma_0(\mathbf p)<\bar\sigma_1(\mathbf p)\leq\bar\sigma_2(\mathbf p)\leq\cdots\). As the iterates approach a target, the first \(m\) positive eigenvalues coalesce and the gradient of the ordered eigenvalue \(\bar\sigma_1\) changes rapidly across the near-multiple set. We monitor the relative spread of the cluster,
\begin{equation}\label{eq:cluster-spread}
\delta_m(\mathbf p):=\frac{\max_{k\in[m]}\bar\sigma_k(\mathbf p)-\min_{k\in[m]}\bar\sigma_k(\mathbf p)}{\tfrac1m\sum_{k\in[m]}\bar\sigma_k(\mathbf p)},
\end{equation}
and switch to the second stage once \(\delta_m\) falls below a prescribed tolerance; the two stages are run once each, in succession, and the second is run to convergence.

\emph{Stage 2: cluster refinement.} Near the active manifold the problem is locally the smooth program
\begin{equation}\label{eq:local-sqp}
\max_{\mathbf p,\,t}\ t\qquad\text{subject to}\qquad \bar{\mathsf S}(\mathbf p)-tI=0,
\end{equation}
with \(\tfrac12m(m+1)\) constraints. The second stage alternates two kinds of steps, both built on the derivatives \eqref{eq:cluster-derivative}, in the pattern of \cite[Algorithm~4.1]{Bareilles_2023}. While the cluster is not yet coalesced, a proximal maximin step in the manner of Overton \cite{Overton_1988} maximizes the smallest eigenvalue of the linearized cluster,
\begin{equation}\label{eq:maximin-subproblem}
\max_{d}\ \lambda_{\min}\Bigl(C+\sum_{j\in[P]}d_j\,\partial_{p_j}\bar{\mathsf S}(\mathbf p)\Bigr)-\frac{\|d\|^2}{2\eta},\qquad C:=\bar{\mathsf S}(\mathbf p)-\tfrac1m\operatorname{tr}\bar{\mathsf S}(\mathbf p)\,I,
\end{equation}
with \(\eta>0\) a proximal parameter; once it has coalesced, a sequential quadratic programming (SQP) step for \eqref{eq:local-sqp} solves
\begin{equation}\label{eq:sqp-subproblem}
\max_{d,\,t}\ t-\tfrac12 d^{\mathsf T}Bd\qquad\text{subject to}\qquad \bar{\mathsf S}(\mathbf p)+\sum_{j\in[P]}d_j\,\partial_{p_j}\bar{\mathsf S}(\mathbf p)=tI,
\end{equation}
with \(B\succ0\) an approximation of the negative Hessian of the Lagrangian. With a symmetric \(m\times m\) Lagrange multiplier \(\Lambda\) the Lagrangian of \eqref{eq:local-sqp} is \(t+\langle\Lambda,\bar{\mathsf S}(\mathbf p)-tI\rangle_F\), where \(\langle A,B\rangle_F:=\operatorname{tr}(A^{\mathsf T}B)\) is the Frobenius inner product, so its stationarity conditions are \(\operatorname{tr}\Lambda=1\) and
\begin{equation}\label{eq:cluster-stationarity}
\bigl\langle\Lambda,\partial_{p_j}\bar{\mathsf S}(\mathbf p^\star)\bigr\rangle_F=0\qquad\text{for every component }p_j\text{ of }\mathbf p .
\end{equation}
The multiplier of either subproblem is the current estimate of \(\Lambda\), and the residual of \eqref{eq:cluster-stationarity},
\begin{equation}\label{eq:kkt-residual}
\mathcal K(\mathbf p,\Lambda):=\max_{j\in[P]}\,\bigl|\bigl\langle\Lambda,\partial_{p_j}\bar{\mathsf S}(\mathbf p)\bigr\rangle_F\bigr|,
\end{equation}
the maximum being over all components of \(\mathbf p=(\tau_1,\tau_2,\mathbf c,\mathbf r,\boldsymbol\alpha,\boldsymbol\beta)\), is the natural measure of convergence: Stage~2 terminates when both \(\delta_m\) and \(\mathcal K\) fall below a prescribed tolerance. Strict complementarity, \(\Lambda\succ0\), is the condition under which Newton steps on \eqref{eq:local-sqp} converge locally quadratically \cite[Theorem~4.5]{Bareilles_2023}. Unlike \cite{Bareilles_2023}, whose Newton step uses the exact Hessian of the Lagrangian, the SQP step is run as a quasi-Newton method on the first derivatives \eqref{eq:cluster-derivative} alone: second derivatives of \(\bar{\mathsf S}\) would require second derivatives of the torus Green function in the modulus and the hole geometry. The switching rule, the solution of the two subproblems, and the step control are given in Section~\ref{ssec:optimizer-settings}.

\subsection{Stationarity conditions}
\label{ssec:stationarity}

The optimization returns a point \(\mathbf p^\star\) satisfying \eqref{eq:intro-three-dimensional-cluster} and a certificate of stationarity \(\Lambda\), a multiplier satisfying \eqref{eq:cluster-stationarity}, both to the prescribed tolerance. Such a point is a Clarke stationary point of \eqref{eq:objective} exactly when \eqref{eq:cluster-stationarity} holds for some \(\Lambda\succeq0\) with \(\operatorname{tr}\Lambda=1\), that is, when \(0\) lies in the Clarke subdifferential of \(\lambda_{\min}\circ\bar{\mathsf S}\) \cite{Lewis_Overton_1996}. The positive semidefiniteness of the multiplier makes the point stationary for the nonsmooth problem and not merely for \eqref{eq:local-sqp}; this is the first-order necessary condition for a local maximizer. The certificate is more than a diagnostic: it contains the rescaling of the eigenfunctions that produces the surface, and the criteria of Section~\ref{ssec:verification} test exactly the conclusions drawn here.

\emph{Output.} We use the superscript \({}^\star\) for the computed candidate: \(\mathbf p^\star=(\tau^\star,\mathbf c^\star,\mathbf r^\star,\boldsymbol\alpha^\star,\boldsymbol\beta^\star)\), \(\rho^\star:=\rho_{\mathbf p^\star}\), \(\Omega^\star:=\Omega_{\mathbf p^\star}\), \(\bar\sigma_1^\star:=\bar\sigma_1(\mathbf p^\star)\), the certificate \(\Lambda^\star\), and \(X^\star:=(L_{\rho^\star}\Lambda^\star)^{1/2}(u_1,\dots,u_m)^{\mathsf T}\), whose role is explained by the following theorem.

\begin{theorem}[Stationarity certificate]\label{prop:certificate}
Let \(\mathbf p^\star\) be a Clarke-stationary point of \eqref{eq:objective} at which \(\sigma_1\) has multiplicity \(m\), let \(\Lambda\) be a certificate in \eqref{eq:cluster-stationarity}, and let \(u_1,\dots,u_m\) be eigenfunctions normalized by \eqref{eq:weighted-eigenfunction-normalization}. Set \(X:=(L_\rho\Lambda)^{1/2}(u_1,\dots,u_m)^{\mathsf T}\colon\overline{\Omega^\star}\to\mathbb R^m\). Then:
\begin{enumerate}
\item[(a)] the components of \eqref{eq:cluster-stationarity} in the weight directions state that \(|X|^2-1\) is orthogonal to every trigonometric polynomial of degree \(\leq N_\rho\) on each \(\Gamma_\ell\) in \(L^2(\partial\Omega,\rho\,ds_\tau)\); in particular, if \eqref{eq:cluster-stationarity} holds for every smooth \(\delta\omega\), then \(|X|=1\) on \(\partial\Omega\);
\item[(b)] if \(|X|=1\) on \(\partial\Omega\), then the components in the directions \((\tau,\mathbf c,\mathbf r)\) are equivalent to the conformality of \(X\), so that \(X\) is a (possibly branched) free boundary minimal immersion of \(\Omega^\star\) into \(\mathbb B^m\), for a metric conformal to \(|dz|^2\) with boundary length element \(\sigma\rho\,ds_\tau\), with \(\operatorname{rank}\Lambda\) the dimension of the smallest ball containing its image.
\end{enumerate}
\end{theorem}

\begin{proof}
By \eqref{eq:cluster-derivative} for the weight, the components of \eqref{eq:cluster-stationarity} in the directions \(\delta\omega\) read
\[\int_{\partial\Omega}(1-L_\rho\,u^{\mathsf T}\Lambda u)\,\rho\,\delta\omega\,ds_\tau=0;\]
for \(\delta\omega\) ranging over the parametrization \eqref{eq:fourier-boundary-weight} this is (a), since \(L_\rho u^{\mathsf T}\Lambda u=|X|^2\); the constant on \(\Gamma_1\), excluded by \(\alpha_{1,0}=0\), is covered by \(\delta\omega\equiv1\), for which the integral equals \(L_\rho(1-\operatorname{tr}\Lambda)=0\) by the orthonormality \eqref{eq:weighted-eigenfunction-normalization}. For (b), \(X\) is harmonic with \(\partial_{\nu_\tau}X=\sigma\rho X\), so under \(|X|=1\) it is a free boundary harmonic map into \(\mathbb B^m\) for any metric \(g\) conformal to \(|dz|^2\) that equals \(\sigma^2\rho^2|dz|^2\) along \(\partial\Omega\) (harmonicity is conformally invariant in two dimensions), and it remains to show that the geometric components of \eqref{eq:cluster-stationarity} force it to be conformal. Let \(\mathcal T(X):=dX^{\mathsf T}\!\otimes dX-\tfrac12|dX|_g^2\,g\) be the trace-free stress--energy tensor of \(X\). For a variation \(h\) of the metric, the first variation formula of \textcite[Lemma~5.1]{Fraser_2015}, in its polarized form for the cluster matrix (compare \eqref{eq:cluster-derivative}), gives \(\partial_h\bar{\mathsf S}_{kk'}=-L_\rho\int_\Omega\langle\mathcal T(u_k,u_{k'}),h\rangle\,dA-\tfrac12\sigma L_\rho\int_{\partial\Omega}u_ku_{k'}h(T,T)\,\rho\,ds_\tau+\sigma\,\partial_hL_\rho\,\delta_{kk'}\), where \(\mathcal T(u,v)\) is the polarization of \(\mathcal T\) and \(T\) the unit tangent. Pairing with \(\Lambda\) and using \(L_\rho u^{\mathsf T}\Lambda u=|X|^2=1\) on \(\partial\Omega\), the two boundary terms cancel, because \(\partial_hL_\rho=\tfrac12\int_{\partial\Omega}h(T,T)\rho\,ds_\tau\), and the geometric components of \eqref{eq:cluster-stationarity} become \(\int_\Omega\langle\mathcal T(X),h\rangle\,dA=0\) for every \(h\) tangent to the family \((\tau,\mathbf c,\mathbf r)\). Since \(X\) is harmonic, \(\mathcal T(X)\) is divergence free, and on \(\partial\Omega\) it satisfies \(\mathcal T(X)(\nu,T)=\langle\partial_\nu X,\partial_TX\rangle=\tfrac12\sigma\rho\,\partial_T|X|^2=0\); hence \(\int_\Omega\langle\mathcal T(X),\mathcal L_Vg\rangle\,dA=0\) for every vector field \(V\) tangent to \(\partial\Omega\), and the pairing depends only on the class of \(h\) in the tangent space of the moduli space. By the dimension count of Section~\ref{ssec:circle-domain}, the directions \((\tau,\mathbf c,\mathbf r)\) with \(c_1\) fixed span this \(3b\)-dimensional space, and since \(\mathcal T(X)\) is trace free the conformal part of \(h\) does not contribute either, so the identity holds for every trace-free \(h\); taking \(h=\mathcal T(X)\) gives \(\mathcal T(X)=0\), i.e.,\ \(X\) is conformal. A conformal free boundary harmonic map is a (branched) free boundary minimal immersion, and \(\operatorname{rank}\Lambda=\operatorname{rank}(L_\rho\Lambda)^{1/2}\) is the dimension of the span of the coordinates of \(X\). This is the argument of \cite[Proposition~5.2]{Fraser_2015} with the convex-hull step of their Lemma~5.3, which produces the weights on the eigenfunctions from maximality, replaced by the multiplier \(\Lambda\). The continuous counterparts of (a) and (b) are \cite[Propositions~2.8 and~2.4]{Fraser_2013} for maximizers in a conformal class and over all metrics, respectively, and \cite[Theorem~3]{Karpukhin_Metras_2022} for conformally extremal pairs \((g,\rho)\) on surfaces, the setting of the weighted problem \eqref{eq:weighted-steklov-problem}.
\end{proof}
Theorem~\ref{prop:certificate} is the finite-dimensional form of the Fraser--Schoen correspondence: the multiplier that certifies stationarity of the nonsmooth problem already contains the rescaling of the eigenfunctions that produces the surface, and the parametrization of Section~\ref{ssec:circle-domain} loses nothing because it is complete. In \cite[\S4.3]{Kao_2020} the matrix that combines the eigenfunctions into the surface is determined after the optimization, by a second, nonconvex least-squares minimization of the sphere and conformality conditions, which reached the sphere condition to \(10^{-3}\). Here no second step is needed. The dual variable of the maximin step \eqref{eq:maximin-subproblem}, equivalently the multiplier of \eqref{eq:sqp-subproblem}, is formed at every iteration and is the surface, and the sphere condition holds to the accuracy of the stationarity residual (Table~\ref{tab:sanity-checks}). It also determines the eigenbasis: \(X\) is fixed up to an orthogonal transformation of \(\mathbb R^m\), a rotation of the surface. Parts (a) and (b) also identify the weight. On \(\partial\Omega\), \(\partial_{\nu_\tau}X=\sigma\rho X\) and \(|X|=1\), so \(|\partial_{\nu_\tau}X|=\sigma\rho\), and by conformality \(|\partial_TX|=\sigma\rho\) as well: the boundary length element of the surface is \(\sigma\rho\,ds_\tau\) (compare \cite{Karpukhin_Metras_2022}). The boundary of the surface therefore has total length
\begin{equation}\label{eq:boundary-length}
L_\partial:=\int_{\partial\Omega}|\partial_TX|\,ds_\tau=\int_{\partial\Omega}\sigma\rho\,ds_\tau=\bar\sigma_1,
\end{equation}
and the Fraser--Schoen bound \(\bar\sigma_1\leq2\pi(\gamma+b)\) \cite{Fraser_2015} bounds the length of the boundary. Finally, \(\lambda_{\min}(\Lambda)>0\) says that all \(m\) coordinates of \(X\) are present, and it is the strict complementarity condition \(\Lambda\succ0\) under which the Newton step on \eqref{eq:local-sqp} converges quadratically (Section~\ref{ssec:two-stage-optimization}).

\subsection{Verification criteria}
\label{ssec:verification}

A map \(X\colon\overline\Omega\to\mathbb R^m\) is a free boundary minimal immersion into \(\mathbb B^m\) if it is (i) an immersion, (ii) conformal, (iii) harmonic, (iv) maps \(\partial\Omega\) into the unit sphere, and (v) meets the sphere orthogonally along \(\partial\Omega\). Conditions (ii) and (iii) together imply that \(X\) is minimal, \(H=0\); conditions (iv) and (v) are the free boundary conditions \eqref{eq:intro-coordinate-steklov}. For the computed candidate \(X^\star\), (iii) holds by construction, up to the error of the boundary-integral eigensolver, since each coordinate is a linear combination of Steklov eigenfunctions, and (v) follows from (iv): the Steklov condition gives \(\partial_{\nu_\tau}X^\star=\sigma\rho X^\star\), so the conormal of the surface is parallel to \(X^\star\), the normal of the sphere, wherever \(|X^\star|=1\). The remaining conditions are verified a posteriori through the following pointwise residuals.

Write \(z=x+iy\) for the flat coordinate on \(\Omega^\star\), and let \(E:=|X^\star_x|^2\), \(F:=\langle X^\star_x,X^\star_y\rangle\), \(G:=|X^\star_y|^2\) be the coefficients of the first fundamental form and \(J:=\sqrt{EG-F^2}\) the area element. Condition (i) is \(J>0\) on \(\overline{\Omega^\star}\). For (ii) we use the scale-free conformality defect
\begin{equation}\label{eq:conformality-defect}
\mathcal C_{X^\star}:=\frac{\sqrt{(E-G)^2+4F^2}}{E+G}=\frac{s_1^2-s_2^2}{s_1^2+s_2^2}\in[0,1],
\end{equation}
where \(s_1\geq s_2\) are the singular values of \(dX^\star\); it vanishes exactly at conformal points and equals one where \(dX^\star\) has rank one. Since \(X^\star\) is harmonic, \(f_{X^\star}:=\sum_{j\in[m]}(\partial_zX^\star_j)^2=\tfrac14\bigl((E-G)-2iF\bigr)\) is holomorphic in \(\Omega^\star\), so \(|f_{X^\star}|=\tfrac14(E+G)\,\mathcal C_{X^\star}\) attains its maximum on \(\partial\Omega^\star\): away from branch points, where \(E+G\) is bounded below, the conformality defect in the interior is controlled by its boundary values, which are computable from the traces and the Steklov condition alone. As an independent check of (ii)--(iii) we also evaluate the mean curvature
\begin{equation}\label{eq:mean-curvature}
H:=\frac{G\langle X^\star_{xx},\mathbf n\rangle-2F\langle X^\star_{xy},\mathbf n\rangle+E\langle X^\star_{yy},\mathbf n\rangle}{2(EG-F^2)},\qquad \mathbf n:=\frac{X^\star_x\times X^\star_y}{J},
\end{equation}
written for \(m=3\); for general \(m\) the unit normal is replaced by the projection onto the normal space of the surface, and \(|H|\) by the length of the mean curvature vector. Condition (iv) is measured by \(\mathcal R_{\mathrm{sph}}:=\bigl||X^\star|-1\bigr|\). The identity \eqref{eq:boundary-length} is checked directly: the length \(L_\partial\) of the boundary curve of \(X^\star\), evaluated from the traces by the trapezoidal rule, is compared with \(\bar\sigma_1^\star\); this tests the surface against the eigenvalue rather than against the sphere. The residuals \(\mathcal C_{X^\star}\) and \(\mathcal R_{\mathrm{sph}}\) are dimensionless, \(L_\partial\) is measured in units of the radius of the ball and \(H\) in units of its inverse. In Section~\ref{sec:experiments} we report the maxima of \(\mathcal R_{\mathrm{sph}}\) and \(\mathcal C_{X^\star}\) over the boundary nodes, the difference \(L_\partial-\bar\sigma_1^\star\), and the maximum of \(|H|\) and the minimum of \(J\) over the interior sample points described in Section~\ref{ssec:visualization}. These criteria test for a free boundary minimal \emph{immersion}, which is what the Fraser--Schoen correspondence guarantees.

\section{Discretization and implementation}
\label{sec:implementation}

\subsection{Discretization of the eigenproblem}
\label{ssec:nystrom}

Each circle \(\Gamma_\ell\) carries \(N_\ell\) nodes equispaced in \(\theta\), with \(N_\ell\) even and \(N:=\sum_{\ell\in[b]}N_\ell\) in total, and the boundary integral equation \eqref{eq:weighted-steklov-bie} is discretized by the Nystr\"om method of \cite{kim2025layer}: the discrete unknown is the vector of nodal values of the density \(\varphi\), the weight enters as the diagonal matrix of nodal values of \(\rho\), and the result is a dense \(N\times N\) generalized eigenvalue problem \(\mathsf A_h\varphi=\sigma\mathsf B_h\varphi\). The kernels of \(\mathcal S\) and \(\mathcal D^*\) are split as in \eqref{eq:sqbx-kernel-split} of Appendix~\ref{sec:qbx}, \(G_\tau(w)=-\tfrac{1}{2\pi}\log|w|+R_\tau(w)\), into the planar logarithmic kernel and a remainder that is smooth near \(w=0\). For \(\mathcal S\) the logarithmic part is integrated with the periodic product quadrature of Kress \cite[\S12.3]{Kress_2014} and the smooth remainder, as well as all hole-to-hole blocks, with the periodic trapezoidal rule. For \(\mathcal D^*\) there is no singularity to treat: on a circle the logarithmic part of the kernel \(\kappa_{\mathcal D^*}\) is the constant \(1/(4\pi r_\ell)\) for every pair of points, \(-1/r_\ell\) being the signed curvature of \(\Gamma_\ell\) as part of \(\partial\Omega\), and the smooth remainder vanishes on the diagonal since \(\nabla R_\tau(0)=0\), so the trapezoidal rule is applied to the whole self-interaction block with the diagonal entry \(1/(4\pi r_\ell)\). Both rules are spectrally accurate for smooth periodic data, so the discrete eigenvalues converge spectrally in \(N_\ell\) \cite{kim2025layer}, and our implementation reproduces the Steklov eigenvalues tabulated in \cite[Table~3]{kim2025layer} to all printed digits. The eigenvalue problem is solved with a dense generalized eigensolver for the full spectrum; \(N\leq1024\) in the optimization and \(N\leq3200\) in the refinement of Section~\ref{ssec:optimizer-settings}, so the cost is small.

The right and left eigenvectors give the nodal values of \(\varphi_k\) and of the trace of \(u_k=\mathcal S_0[\varphi_k]\), the quadrature weights \(w_i=\rho_i\,ds_i\) discretize the inner product \eqref{eq:weighted-eigenfunction-normalization}, and tangential derivatives of the traces are obtained spectrally from their Fourier coefficients. For a cluster of \(m\) eigenvalues the right eigenvectors \(\Phi\) are orthonormalized in this inner product through their traces, \(\Phi\leftarrow\Phi\,\mathsf G^{-1/2}\) with \(\mathsf G\) the weighted Gram matrix of the traces, and the left eigenvectors are replaced by the dual basis \(\Psi^{\mathsf T}\leftarrow(\Psi^{\mathsf T}\mathsf B_h\Phi)^{-1}\Psi^{\mathsf T}\), which gives the normalization \(\Psi^{\mathsf T}\mathsf B_h\Phi=I\) of Section~\ref{ssec:two-stage-optimization}.

The derivative matrices \(\partial_p\mathsf A_h\), \(\partial_p\mathsf B_h\) of Proposition~\ref{prop:unified-derivative} are assembled from the same node set with the kernel derivatives in place of the kernels and with the same quadrature rules, since the singular parts of the differentiated kernels are of the same type as those of the kernels themselves. On the self-interaction block of \(\Gamma_\ell\), where \(w=r_\ell(e^{i\theta}-e^{i\theta'})\), the logarithmic part of \(G_\tau\) enters \(\partial_{r_\ell}\mathcal S\) and \(\partial_{r_\ell}\mathcal D^*\) only through \(w\cdot\nabla_w(-\tfrac1{2\pi}\log|w|)=-\tfrac1{2\pi}\) and the \(r_\ell\)-independent weight \(\tfrac1{4\pi}\,d\theta'\), and the diagonal entries of the first two terms of \eqref{eq:radius-operator-derivative-K} are \(-1/(2\pi r_\ell)\) for \(\mathcal S\) and \(-1/(4\pi r_\ell^2)\) for \(\mathcal D^*\), the smooth remainder \(R_\tau\) contributing nothing there. The Hessian of \(G_\tau\), needed for \(\partial_{c_\ell}\mathcal D^*\) and \(\partial_{r_\ell}\mathcal D^*\), is \(-\tfrac1{2\pi}\) times the Hessian of \(\log|\vartheta_1|\) plus the constant Hessian of \((\Im z)^2/(2\Im\tau)\), from the theta series or the Weierstrass \(\wp\)-function (Appendix~\ref{sec:analytic formula for dtau}). The modulus derivative kernel \(\partial_\tau G_\tau\) is evaluated from the theta series of Appendix~\ref{sec:analytic formula for dtau}, truncated when the relative contribution of a term falls below \(10^{-14}\); the removable diagonal limit of \(\partial_\tau\vartheta_1/\vartheta_1\) is evaluated as \(\partial_{\tau z}\vartheta_1/\partial_z\vartheta_1\), and since \(\partial_\tau G_\tau\) has no logarithmic singularity, \(\dot{\mathcal S}\) and \(\dot{\mathcal D}^*\) are assembled with the trapezoidal rule alone.

\subsection{Optimization implementation}
\label{ssec:optimizer-settings}

\emph{Stage 1.} The gradient ascent of Section~\ref{ssec:two-stage-optimization} uses a decaying learning rate and feasibility backtracking. The trial step from the iterate \(\mathbf p^{(s)}\) at backtracking level \(\mu\) is \(\mathbf p_{\mathrm{trial}}=\mathbf p^{(s)}+\tfrac{\xi_{\mathrm{lr}}\,2^{-\mu}}{1+\gamma_{\mathrm{lr}}s}\nabla\bar\sigma_1(\mathbf p^{(s)})\), \(\mu=0,\dots,8\), where \(\xi_{\mathrm{lr}}=10^{-2}\) is the initial learning rate and \(\gamma_{\mathrm{lr}}=0.03\) is the decay. A trial step is accepted if it satisfies the constraints of \eqref{eq:objective}, with a small safety margin that keeps the quadrature accurate near the boundary of the admissible set, and the eigenproblem at \(\mathbf p_{\mathrm{trial}}\) can be solved; the margins are never active at the configurations of Section~\ref{sec:experiments}. The longest feasible of these steps is taken. Stage~1 terminates when the relative cluster spread \(\delta_m\) in \eqref{eq:cluster-spread} falls below \(10^{-4}\).

\emph{Stage 2.} The maximin step \eqref{eq:maximin-subproblem} is taken while
\begin{equation}\label{eq:cluster_spread_threshold}
\max_{k\in[m]}\bar\sigma_k(\mathbf p^{(s)})-\min_{k\in[m]}\bar\sigma_k(\mathbf p^{(s)})>\epsilon_{\mathrm{rel}}\,\bar\sigma_1(\mathbf p^{(s)}),\qquad \epsilon_{\mathrm{rel}}:=10^{-9},
\end{equation}
and the SQP step \eqref{eq:sqp-subproblem} once \eqref{eq:cluster_spread_threshold} fails. Both use the backtracking \(\mathbf p_{\mathrm{trial}}=\mathbf p^{(s)}+2^{-\mu}d\), \(\mu=0,\dots,8\), along the computed direction \(d\): the longest trial that is feasible in the sense of Stage~1 and does not decrease \(\min_{k\in[m]}\bar\sigma_k\) is taken, the second test being the fallback rule of \cite[Algorithm~4.1]{Bareilles_2023}, and if every level fails \(\mathbf p^{(s)}\) is retained. Stage~2 stops when \(\delta_m\) and \(\mathcal K\) fall below \(10^{-9}\).

\emph{Proximal maximin step.} We solve \eqref{eq:maximin-subproblem} through its dual
\begin{equation}\label{eq:dual}
\underset{\Lambda\succeq0,\ \operatorname{tr}\Lambda=1}{\operatorname{minimize}}\quad\langle\Lambda,C^{(s)}\rangle_F+\frac{\eta^{(s)}}{2}\sum_{j\in[P]}\langle\Lambda,\partial_{p_j}\bar{\mathsf S}(\mathbf p^{(s)})\rangle_F^2,
\end{equation}
a convex quadratic program over the \(m\times m\) positive semidefinite matrices of unit trace, by projected gradient descent, the projection mapping the eigenvalues of the trial matrix onto the unit simplex; the duality gap of \eqref{eq:dual} is checked, and the direction is \(d_j=\eta^{(s)}\langle\Lambda^{(s)},\partial_{p_j}\bar{\mathsf S}(\mathbf p^{(s)})\rangle_F\). A trial is required in addition to achieve \(5\%\) of the predicted improvement \(\Delta_{\mathrm{pred}}:=\lambda_{\min}\bigl(C^{(s)}+2^{-\mu}\sum_{j\in[P]}d_j\,\partial_{p_j}\bar{\mathsf S}(\mathbf p^{(s)})\bigr)-\lambda_{\min}(C^{(s)})\). The proximal parameter \(\eta^{(s)}\) is adapted by a standard trust-region rule from the ratio of the actual gain \(\bar\sigma_1(\mathbf p_{\mathrm{trial}})-\bar\sigma_1(\mathbf p^{(s)})\) to \(\Delta_{\mathrm{pred}}\).

\emph{SQP step.} Where \textcite{Bareilles_2023} use the exact Hessian of the Lagrangian in \eqref{eq:sqp-subproblem}, we use a positive definite approximation \(B_s\approx-\bigl[\langle\Lambda^{(s)},\partial_{p_i}\partial_{p_j}\bar{\mathsf S}(\mathbf p^{(s)})\rangle_F\bigr]_{ij}\) built by Powell-damped BFGS updates \cite[Procedure~18.2]{Nocedal_Wright_2006} from the displacements \(\mathbf p^{(s+1)}-\mathbf p^{(s)}\) and the corresponding differences of \(-\langle\Lambda^{(s)},\partial_{p_j}\bar{\mathsf S}\rangle_F\), initialized at \(\max\{1/\eta^{(s)},10^{-6}\}I\) when the SQP phase is first entered; no second derivatives of \(\bar{\mathsf S}\) are formed. Since the coalescence manifold is curved, a full step satisfies only \(\bar{\mathsf S}(\mathbf p^{(s)}+d)=tI+O(\|d\|^2)\); if the recomputed spread exceeds \(10^{-10}\), the second-order correction of \cite[\S5.2]{Bareilles_2023}, the minimum-norm \(d_{\mathrm{corr}}\) solving the linearized constraint at the trial point with derivatives evaluated there, is applied before backtracking, subject to the usual trust-region and non-deterioration safeguards; after a step at which every backtracking level failed, the Hessian approximation is regularized and the trust radius reduced. The result is a safeguarded quasi-Newton refinement on the active manifold; the local quadratic convergence proved in \cite[Theorem~4.5]{Bareilles_2023} for the exact Newton step is not claimed for it.

\emph{Refinement.} In a symmetry class such as that of Section~\ref{ssec:setup} the normalized eigenvalues of the cluster form \(q\leq m\) distinct branches \(v_1,\dots,v_q\), one per irreducible representation, each a smooth function of \(\mathbf p\) as long as the symmetry is kept, so that stationarity on the coalescence manifold becomes the smooth system \(v_a(\mathbf p)-v_1(\mathbf p)=0\) for \(a=2,\dots,q\) and \(\sum_{a\in[q]}\lambda_a\,\partial_{\mathbf p}v_a(\mathbf p)=0\) with \(\lambda_a\geq0\), \(\sum_{a\in[q]}\lambda_a=1\), the \(\lambda_a\) being the eigenvalues of \(\Lambda\) in the symmetry-adapted basis. Once Stage~2 has stopped we solve this system for \((\mathbf p,\lambda_1,\dots,\lambda_q)\) by Gauss--Newton, starting from the Stage~2 output, with the resolution raised to \(N_\ell=400\) and \(N_\rho\in\{56,48,56,64,80,120\}\) for \(b\in\{1,2,3,4,5,8\}\), the added Fourier coefficients initialized at zero; the branches are identified from the symmetry types of the eigenfunctions. The step requires the symmetry, which the two-stage method does not; it reduces the stationarity residual of the reported configurations from \(10^{-9}\) to rounding.

\subsection{Verification and visualization}
\label{ssec:visualization}

The boundary residuals of Section~\ref{ssec:verification} are evaluated at the quadrature nodes from boundary data: \(\mathcal R_{\mathrm{sph}}\) from the traces and \(\mathcal C_{X^\star}\) from the tangential derivatives of the traces and \(\partial_{\nu_\tau}X^\star=\sigma\rho X^\star\).  By the maximum principle for \(f_{X^\star}\), the boundary maximum of \(|f_{X^\star}|=\tfrac14(E+G)\mathcal C_{X^\star}\) bounds its interior values, so that \(\mathcal C_{X^\star}\) in the interior is controlled by its boundary maximum and \(\min_{\Omega^\star}(E+G)\); the boundary maximum of \(\mathcal C_{X^\star}\) is the value reported.

Interior values of \(u_k\) and \(\nabla u_k\), needed for \(J\) and \(H\) and for the plots, are computed from the layer-potential representation \eqref{eq:eigenfunction-layer-representation} on a cell-centered \(120\times120\) grid in the fundamental parallelogram, restricted to \(\Omega^\star\). At targets within \(2.5\,h_\ell\) of a periodic copy of the \(\ell\)-th circle, where \(h_\ell:=2\pi r_\ell^\star/N_\ell\) is the node spacing, the plain Nystr\"om quadrature loses accuracy and the contribution of that component is instead evaluated by quadrature by expansion (Appendix~\ref{sec:qbx}). First and second derivatives of the layer potentials are evaluated analytically, by differentiating \eqref{eq:eigenfunction-layer-representation} under the integral sign with the gradient and Hessian of \(G_\tau\) (Appendix~\ref{sec:analytic formula for dtau}) for the plain quadrature and by differentiating the expansion \eqref{eq:sqbx-exact-expansion} for close targets, so that \(H\) carries no finite-difference error.

For the figures, the surface is \(X^\star=R(L_{\rho^\star}\Lambda^\star)^{1/2}u\) with \(R\in O(3)\) chosen so that the symmetry axes of the surface coincide with the coordinate axes, the axis of the central neck (the common normal of the two sheets for \(b\geq2\)) being the vertical \(x_3\)-axis; each surface in Figure~\ref{fig:surfaces} is shown from an oblique viewpoint, along the \(x_3\)-axis, and from the side. The three components of \(X^\star\) are the coordinate functions plotted in Figure~\ref{fig:eigenfunction_contours}. \(X^\star\) is evaluated at the boundary nodes, on the interior grid, and along the four edges of the parallelogram, the periodic seams of the torus, which are drawn to show how the surface closes up in \(\mathbb R^3\).

\section{Computational experiments and results}
\label{sec:experiments}

\subsection{Setup and initialization}
\label{ssec:setup}

The computations of Sections~\ref{ssec:extremal-configurations}--\ref{ssec:surfaces} are carried out in a symmetry class: the modulus is purely imaginary, so the torus is rectangular, and the configuration is invariant under the reflection in a closed geodesic through the hole centers, the reflection across a perpendicular geodesic through a hole center and, for \(b\geq2\), the translation by \(1/b\) along the geodesic, so that the holes, with equal radii and identical weights, are equispaced on a common closed geodesic. We call these configurations \emph{collinear}, for \(b=1\) as well, where the condition on the center is empty and the symmetry is that of the modulus and the weight. In the implementation the centers are held at \(c_\ell=\tfrac{2\ell-1}{2b}+\tfrac i2\), \(\ell\in[b]\), equispaced on a horizontal closed geodesic, so that \(\bar c_1=\tfrac1{2b}+\tfrac i2\) in \eqref{eq:objective}, the real part of the modulus is held at zero, a single weight profile \(\omega\) with \(\alpha_{1,0}=0\) is shared by all holes, and \(\Im\tau\), the radii, and the remaining Fourier coefficients of \(\omega\) are optimized from \(\tau^{(0)}=i\), \(\omega^{(0)}\equiv0\) and equal radii \(r^{(0)}=0.3\), \(0.15\), \(0.1\), \(0.08\), \(0.06\), \(0.04\) for \(b=1,2,3,4,5,8\). The restriction is lifted for the certificate: the stationarity residual \(\mathcal K\) of Table~\ref{tab:sanity-checks} is evaluated over all components of \(\mathbf p\), the centers and the individual weights included, so the reported configurations are stationary in the full parameter space. Random initializations (random Fourier coefficients of \(\omega\), centers not on a common closed geodesic, \(\tau^{(0)}\) off the imaginary axis) lead to the configurations of Section~\ref{ssec:equilateral}, where the restriction is not imposed.

The two-stage optimization is run with \(N_\ell=128\) nodes per boundary component and \(N_\rho=12\) Fourier modes per component for the weight (Section~\ref{ssec:nystrom}); the configurations reported in Table~\ref{tab:sanity-checks} are the output of the refinement step of Section~\ref{ssec:optimizer-settings}, at \(N_\ell=400\) and \(N_\rho\in\{56,48,56,64,80,120\}\). Increasing each \(N_\rho\) by \(32\), with the stationarity equations re-solved, leaves the printed digits of \(\bar\sigma_1^\star\), \(\Im\tau^\star\), and \(r^\star\) unchanged, and increasing \(N_\ell\) from \(400\) to \(800\) changes \(\bar\sigma_1^\star\) by less than \(10^{-12}\) and keeps \(\mathcal K\) below \(10^{-12}\). The conformality defect is the residual most sensitive to \(N_\rho\): for \(b=8\) it is \(4\times10^{-9}\) at \(N_\rho=80\) and \(10^{-11}\) at \(N_\rho=120\).

\subsection{Extremal collinear configurations}
\label{ssec:extremal-configurations}

\begin{figure}[p]
    \centering
    \includegraphics[width=0.573\textwidth,trim={138 30 80 35},clip]{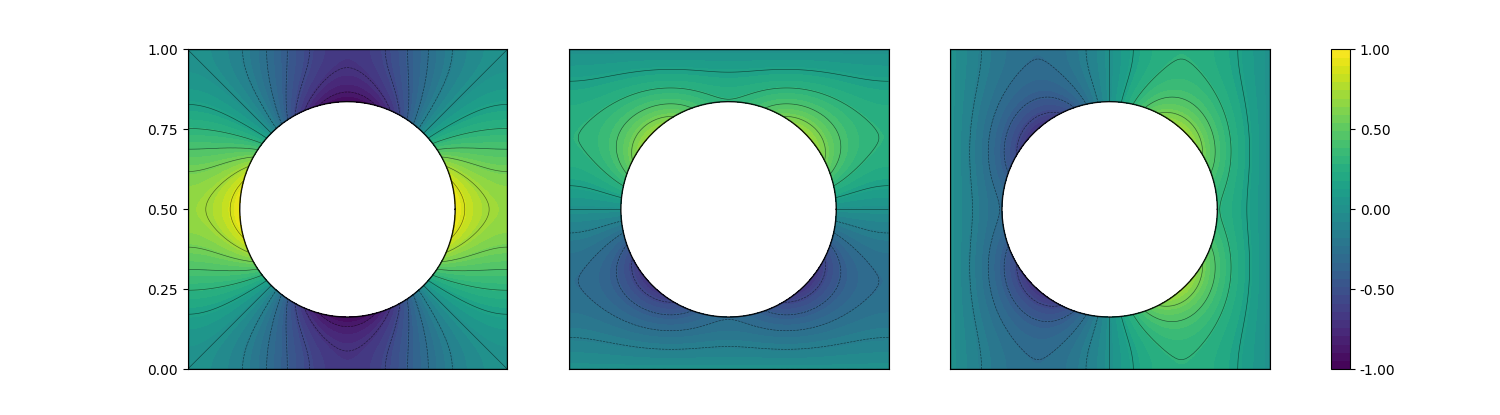}\\
    \includegraphics[width=0.573\textwidth,trim={138 20 80 24},clip]{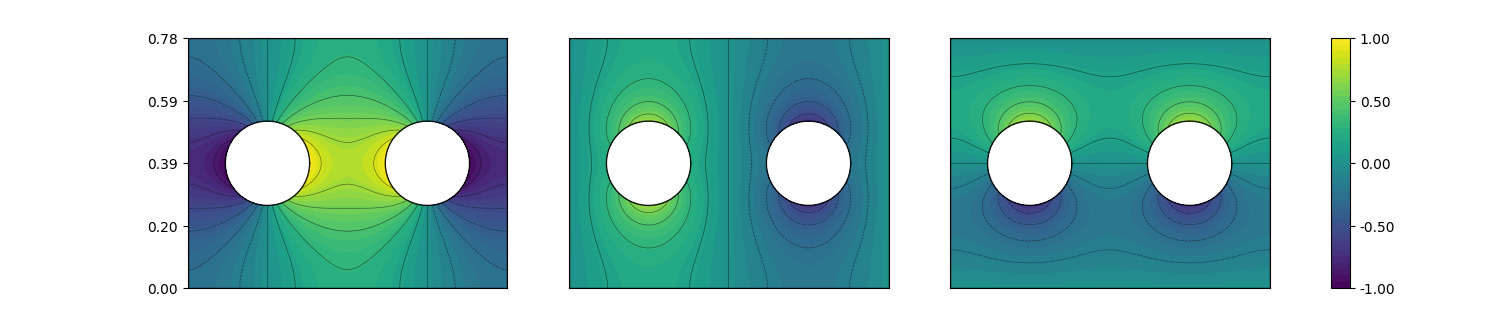}\\
    \includegraphics[width=0.573\textwidth,trim={138 35 80 41},clip]{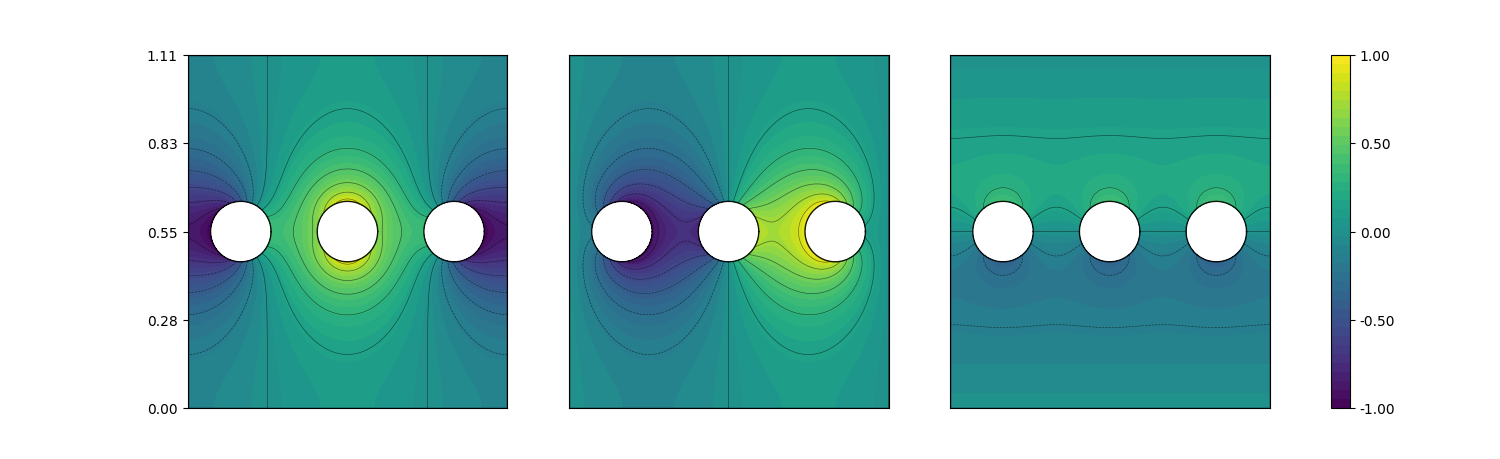}\\
    \includegraphics[width=0.573\textwidth,trim={138 48 80 55},clip]{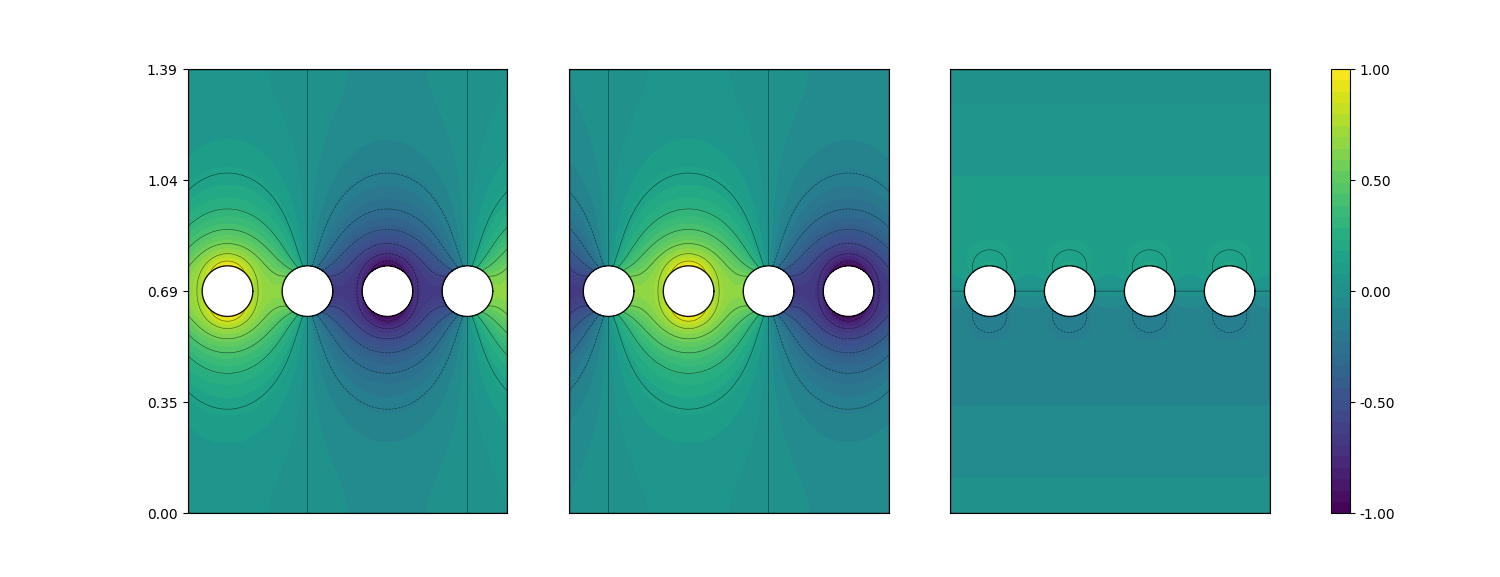}\\
    \includegraphics[width=0.573\textwidth,trim={138 59 80 67},clip]{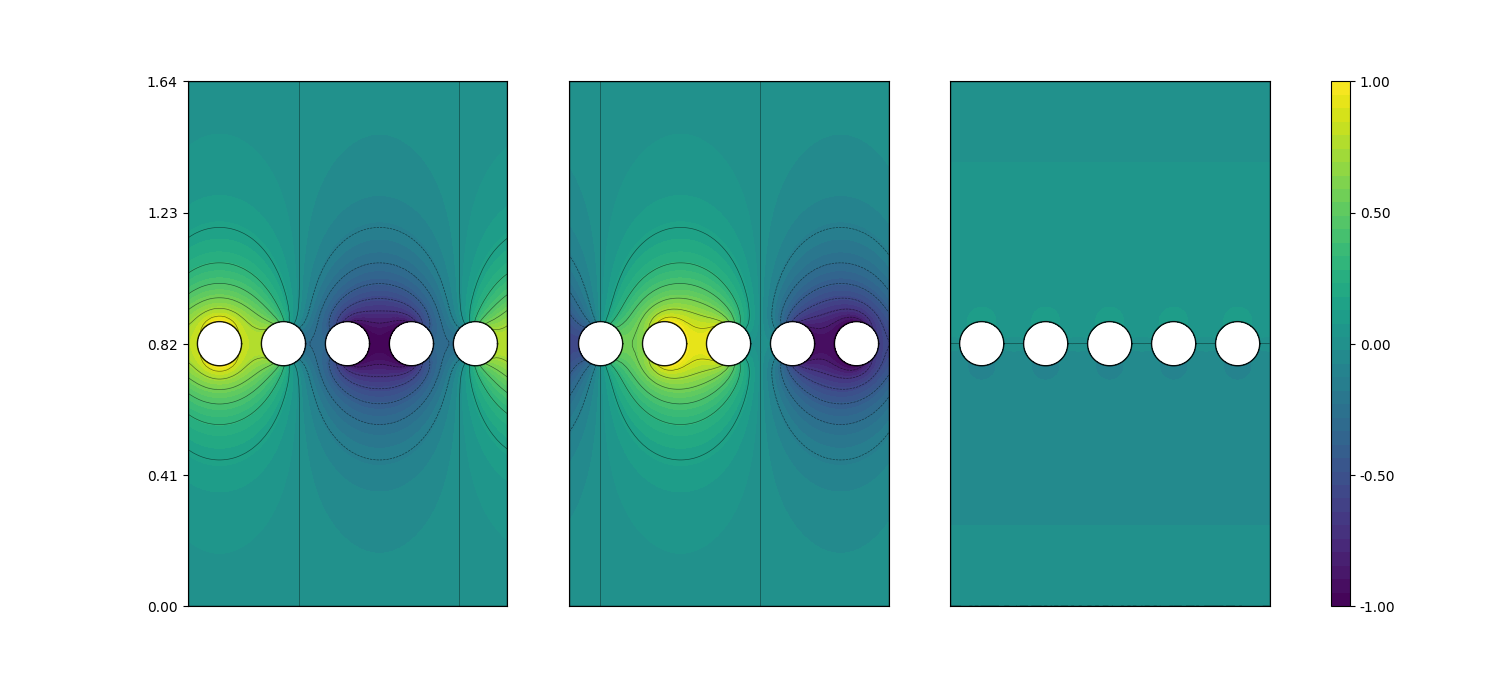}\\
    \includegraphics[width=0.573\textwidth,trim={138 75 80 99},clip]{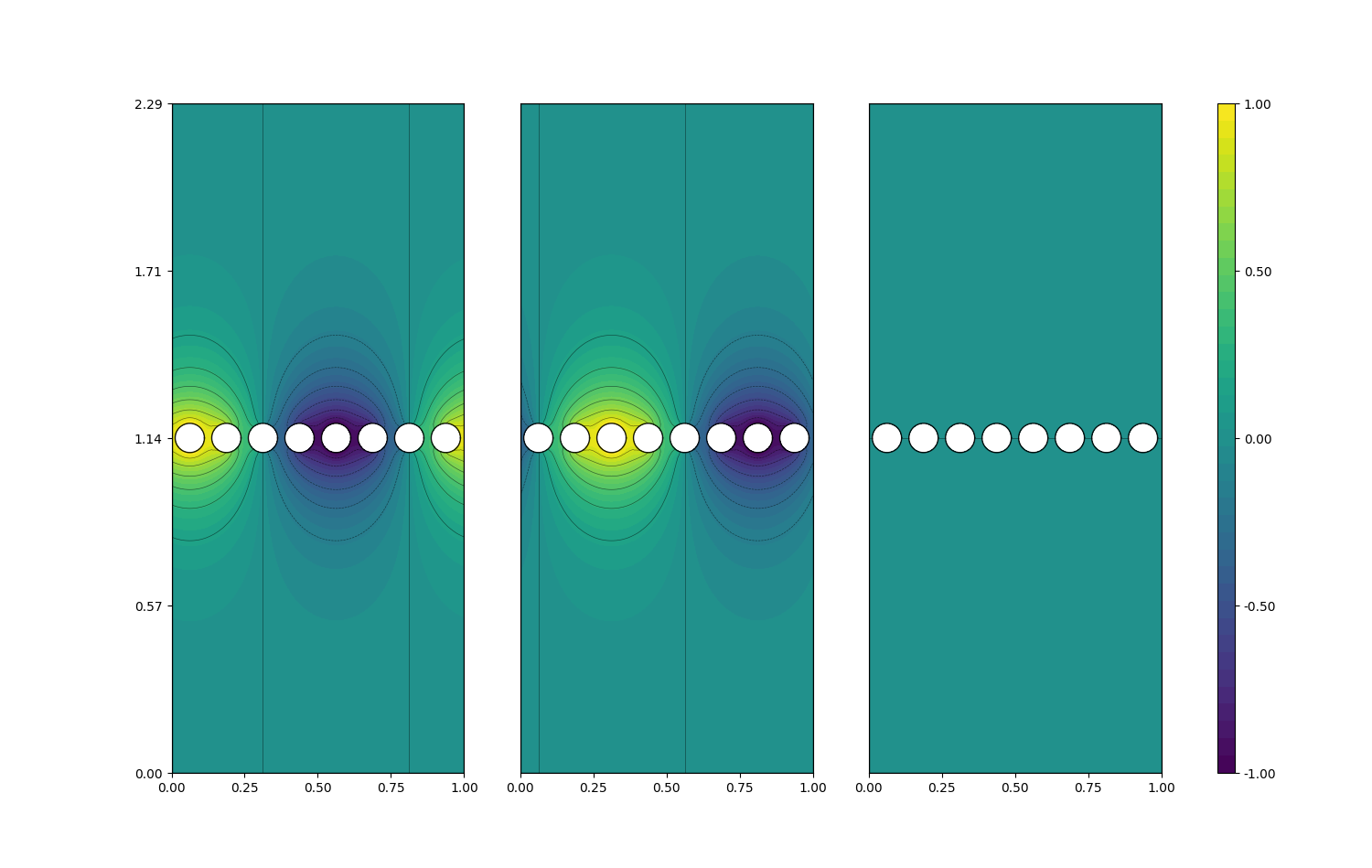}
    \caption{The coordinate functions \(X^\star_1,X^\star_2,X^\star_3\) of the spectral map on \(\overline{\Omega^\star}\) for \(b\in\{1,2,3,4,5,8\}\), one row per \(b\), on the fundamental parallelogram of the computed modulus \(\tau^\star\) with the holes in white (Section~\ref{ssec:extremal-configurations}).}
  \label{fig:eigenfunction_contours}
\end{figure}

\begin{figure}[t]
    \centering
    \includegraphics[width=0.33\textwidth]{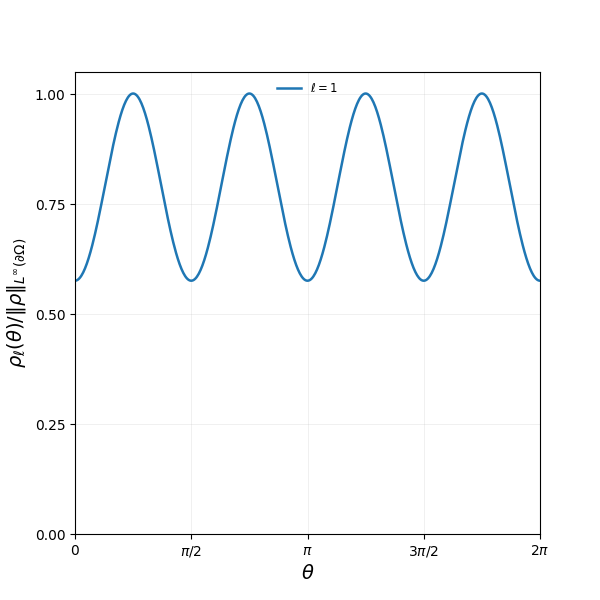}\hfill
    \includegraphics[width=0.33\textwidth]{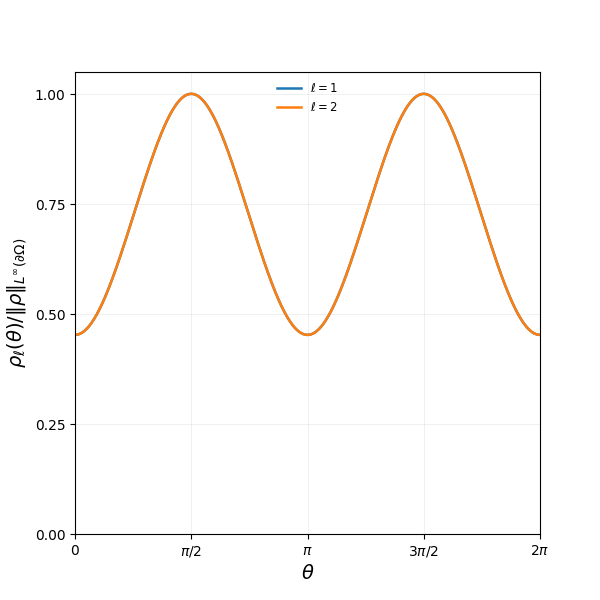}\hfill
    \includegraphics[width=0.33\textwidth]{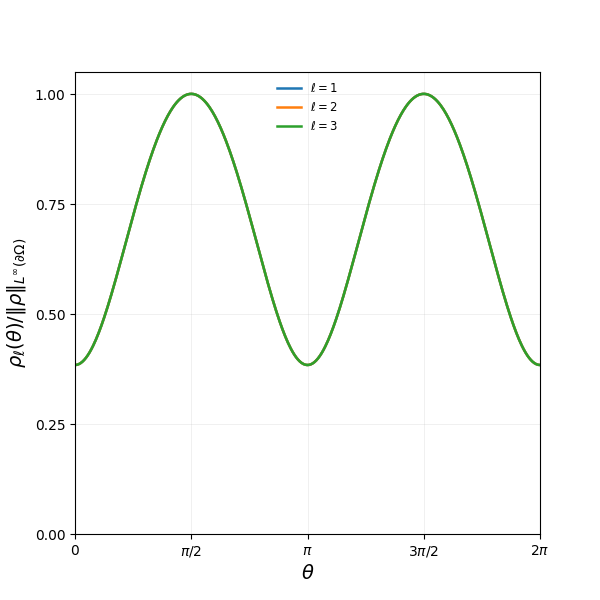}\\
    \includegraphics[width=0.33\textwidth]{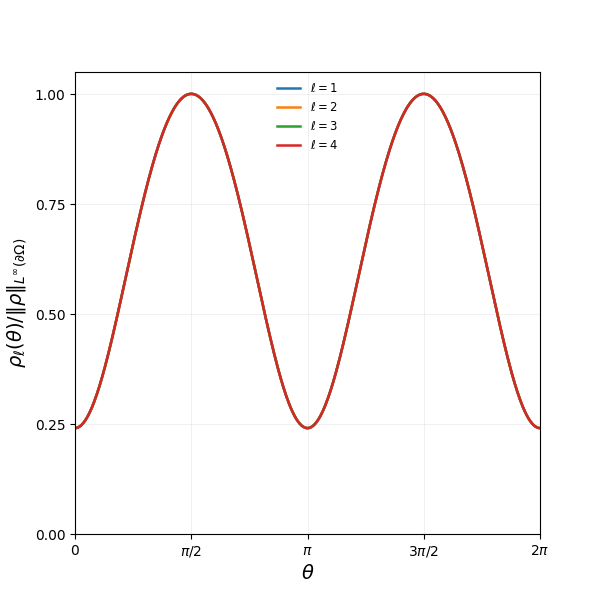}\hfill
    \includegraphics[width=0.33\textwidth]{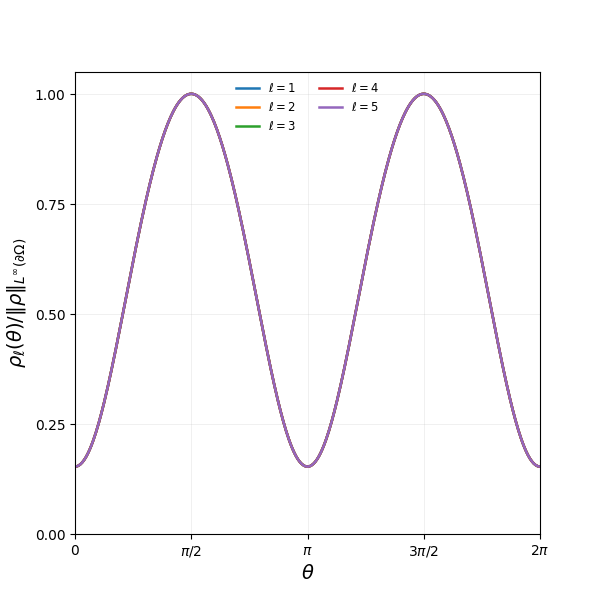}\hfill
    \includegraphics[width=0.33\textwidth]{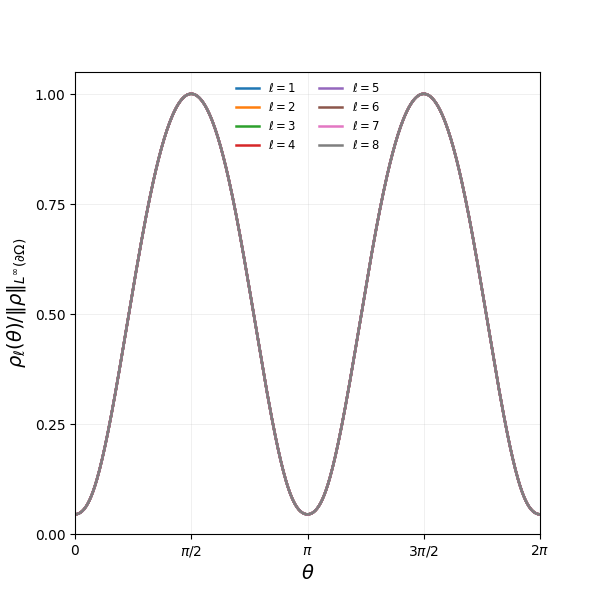}
    \caption{The extremal boundary weight \(\rho^\star(\theta)\), normalized by its maximum, for \(b\in\{1,2,3,4,5,8\}\) (Section~\ref{ssec:extremal-configurations}); the weight is shared by the \(b\) components (Section~\ref{ssec:setup}).}
  \label{fig:boundary_densities}
\end{figure}

At a collinear configuration the first three normalized eigenvalues coalesce with a certificate \(\Lambda^\star\) of full rank, so Theorem~\ref{prop:certificate} produces a map \(X^\star\colon\overline{\Omega^\star}\to\mathbb B^3\); the residuals of Table~\ref{tab:sanity-checks} confirm that \(X^\star\) is a free boundary minimal immersion to the accuracy expected from the discretization, and Figure~\ref{fig:surfaces} shows the surfaces. For \(b=1\) the surface is embedded with connected boundary; for \(b\geq2\) it consists of two sheets joined by a central neck and by \(b\) equatorial necks (Section~\ref{ssec:surfaces}). The top block of Table~\ref{tab:sanity-checks} lists the extremal value \(\bar\sigma_1^\star\) and the modulus \(\tau^\star\), Figure~\ref{fig:eigenfunction_contours} shows the corresponding domains and Figure~\ref{fig:boundary_densities} the weights.

Along the family the modulus is purely imaginary and, for \(b\geq2\), \(\Im\tau^\star\) increases monotonically with \(b\) while the common radius decreases (Table~\ref{tab:sanity-checks}), with \(2br^\star=0.53,\,0.57,\,0.63,\,0.69,\,0.80\) for \(b=2,3,4,5,8\): the torus becomes long in the direction transverse to the geodesic and the holes become small compared to the period. In the surface, a taller fundamental domain and smaller holes correspond to smaller equatorial necks and two sheets that are closer together, which is the regime of the gluing construction of \cite{Folha_Pacard_Zolotareva_2017}. For \(b\geq2\) each weight has period \(\pi\), with maxima at \(\theta=\pm\tfrac\pi2\), the points of \(\Gamma_\ell\) farthest from the geodesic, and minima at \(\theta\in\{0,\pi\}\), the points facing the neighboring holes; the ratio \(\max\rho^\star/\min\rho^\star\) grows from about \(2\) for \(b=2\) to over \(20\) for \(b=8\). The profiles reflect the shape of the surface, since \(\sigma\rho^\star\,ds_\tau\) is the arclength of its boundary (Section~\ref{ssec:stationarity}). Each boundary component consists of two long arcs, one on each sheet, joined by two short passages across the equatorial necks; the arcs are the maxima of the weight and the passages its minima, and the ratio grows with \(b\) because the necks shrink. For \(b\geq2\) the coordinate functions \(X^\star_1,X^\star_2\) are, to leading order, \(\cos(2\pi x)\) and \(\sin(2\pi x)\) modulated along the geodesic, which maps to the equator of \(\mathbb S^2\), while \(X_3^\star\) is odd under the reflection in the geodesic and, for large \(b\), small: the two sheets of the surface are close to the equatorial plane (right column of Figure~\ref{fig:eigenfunction_contours}). The case \(b=1\) is different: \(\tau^\star=i\) to eight digits, the weight has period \(\tfrac\pi2\), and domain and weight are invariant under the full symmetry group of the square torus with one hole, of order eight.

The extremal values \(\bar\sigma_1^\star\) increase with \(b\) (Table~\ref{tab:sanity-checks}). As \(b\to\infty\) the surfaces of \textcite{Folha_Pacard_Zolotareva_2017} converge to the doubly covered unit disk, whose boundary length is \(4\pi\) and on which the coordinate functions have eigenvalue one, consistent with \(\bar\sigma_1^\star\to4\pi\approx12.5664\) along this family, and \(\bar\sigma_1^\star(b{=}8)\) is already within \(3\times10^{-4}\) of the limit.

\begin{table}[t]
\centering
\footnotesize
\setlength{\tabcolsep}{4pt}
\renewcommand{\arraystretch}{0.95}

\sisetup{round-mode=places, round-precision=1, table-format=1.1e+2, group-minimum-digits=20}
\begin{tabular}{@{}l *{6}{S} @{}}
\toprule
& {\(b=1\)} & {\(b=2\)} & {\(b=3\)} & {\(b=4\)} & {\(b=5\)} & {\(b=8\)}\\
\midrule
\multicolumn{7}{@{}l}{\textit{Computed configuration}}\\
\addlinespace[1pt]
\(N_\rho\) & {56} & {48} & {56} & {64} & {80} & {120}\\
\(\bar\sigma_1^\star\) & {\tablenum[scientific-notation=false,round-mode=none,table-format=2.10]{11.3300646633}} & {\tablenum[scientific-notation=false,round-mode=none,table-format=2.10]{11.6611747247}} & {\tablenum[scientific-notation=false,round-mode=none,table-format=2.10]{12.3540550992}} & {\tablenum[scientific-notation=false,round-mode=none,table-format=2.10]{12.5192261877}} & {\tablenum[scientific-notation=false,round-mode=none,table-format=2.10]{12.5550325108}} & {\tablenum[scientific-notation=false,round-mode=none,table-format=2.10]{12.5661315560}} \\
\(\Im\tau^\star\) & {\tablenum[scientific-notation=false,round-mode=none,table-format=2.10]{1.0000000000}} & {\tablenum[scientific-notation=false,round-mode=none,table-format=2.10]{0.7814633148}} & {\tablenum[scientific-notation=false,round-mode=none,table-format=2.10]{1.1058050913}} & {\tablenum[scientific-notation=false,round-mode=none,table-format=2.10]{1.3899779310}} & {\tablenum[scientific-notation=false,round-mode=none,table-format=2.10]{1.6401265921}} & {\tablenum[scientific-notation=false,round-mode=places,round-precision=10,table-format=2.10]{2.2858326420}}\\
\(r^\star\) & {\tablenum[scientific-notation=false,round-mode=none,table-format=2.10]{0.3366882157}} & {\tablenum[scientific-notation=false,round-mode=none,table-format=2.10]{0.1317424031}} & {\tablenum[scientific-notation=false,round-mode=none,table-format=2.10]{0.0944889217}} & {\tablenum[scientific-notation=false,round-mode=none,table-format=2.10]{0.0791200142}} & {\tablenum[scientific-notation=false,round-mode=none,table-format=2.10]{0.0690076567}} & {\tablenum[scientific-notation=false,round-mode=none,table-format=2.10]{0.0499956650}} \\
\midrule
\multicolumn{7}{@{}l}{\textit{Stationarity (Section~\ref{ssec:two-stage-optimization})}}\\
\addlinespace[1pt]
\(\delta_m(\mathbf p^\star)\) \eqref{eq:cluster-spread} & 1.20722591e-14 & 2.34589533e-14 & 2.73195964e-15 & 8.51341838e-15 & 2.12228464e-15 & 4.52354160e-15 \\
\(\mathcal K(\mathbf p^\star,\Lambda^\star)\) \eqref{eq:kkt-residual} & 8.10462808e-14 & 4.02801115e-14 & 4.53526106e-14 & 3.36169109e-14 & 7.65839225e-14 & 8.51159421e-14 \\
\(\lambda_{\min}(\Lambda^\star)\) & 2.8086123077e-01 & 2.6480143946e-01 & 9.0225406463e-02 & 2.4366454859e-02 & 6.7248133675e-03 & 1.950671e-04\\
\midrule
\multicolumn{7}{@{}l}{\textit{Free boundary conditions on \(\partial\Omega^\star\) (Section~\ref{ssec:verification})}}\\
\addlinespace[1pt]
\(\max\mathcal R_{\mathrm{sph}}\) & 1.43218770e-14 & 1.18793864e-14 & 1.67643677e-14 & 5.10702591e-15 & 5.32907052e-15 & 5.55111512e-15 \\
\(L_\partial-\bar\sigma_1^\star\) \eqref{eq:boundary-length} & 3.55271368e-14 & 1.33226763e-13 & 6.03961325e-14 & 3.55271368e-15 & 1.33226763e-13 & 6.89226454e-13 \\
\(\max_{\partial\Omega^\star}\mathcal C_{X^\star}\) \eqref{eq:conformality-defect} & 4.05707108e-13 & 1.23820846e-12 & 1.00079351e-12 & 1.39958703e-12 & 1.62966564e-12 & 9.56770981e-12 \\
\midrule
\multicolumn{7}{@{}l}{\textit{Minimality and immersion in \(\Omega^\star\) (Section~\ref{ssec:verification})}}\\
\addlinespace[1pt]
\(\max|H|\) \eqref{eq:mean-curvature} & 1.24335760e-13 & 1.17444595e-13 & 8.50697985e-13 & 1.50309150e-12 & 6.43421113e-12 & 4.32976228e-09 \\
\(\min J\) & 2.2176524356e+00 & 1.6329948864e+00 & 3.0828290722e-01 & 5.1918375036e-02 & 1.0308913605e-02 & 1.538724e-04\\
\bottomrule
\end{tabular}
\caption{Computed configurations, stationarity measures, and verification residuals for \(b\in\{1,2,3,4,5,8\}\) (Sections~\ref{ssec:extremal-configurations} and~\ref{ssec:verification-results}). The modulus is purely imaginary in every case, so only \(\Im\tau^\star\) is listed; all entries are computed at \(N_\ell=400\) and the listed \(N_\rho\) (Section~\ref{ssec:optimizer-settings}). Boundary residuals are maxima over the quadrature nodes; \(|H|\) and the area element \(J\) are evaluated on the interior grid of Section~\ref{ssec:visualization}. The stationarity residual \(\mathcal K\) is evaluated over the full parameter vector: both components of \(\tau\), all centers and radii, and the Fourier coefficients of \(b\) independent weight profiles through order \(N_\rho\).}
\label{tab:sanity-checks}
\end{table}

\subsection{Verification criteria}
\label{ssec:verification-results}

The stationarity block of Table~\ref{tab:sanity-checks} records how well each configuration satisfies the first-order conditions of Section~\ref{ssec:two-stage-optimization}. For every \(b\) the radii remain equal, which is not imposed, and the stationarity residual \(\mathcal K\), taken over all parameters, shows that the configuration is stationary with respect to the centers and the individual weights as well (Section~\ref{ssec:setup}). The relative cluster spread \(\delta_m\) is below \(10^{-13}\) for every \(b\), so the first three eigenvalues coalesce to rounding accuracy; the stationarity residual \(\mathcal K\) is below \(10^{-13}\) for every \(b\); and \(\lambda_{\min}(\Lambda^\star)\) is positive, so \(\Lambda^\star\) has rank three, all three coordinates of \(X^\star\) are present, and strict complementarity holds.

Since \(\operatorname{tr}\Lambda^\star=1\), the eigenvalues of \(\Lambda^\star\) sum to one; \(\lambda_{\min}(\Lambda^\star)\) decreases monotonically with \(b\), from \(0.28\) for \(b=1\) to \(2.0\times10^{-4}\) for \(b=8\), three orders of magnitude below the other two eigenvalues. The decreasing eigenvalue belongs to the third coordinate: as the two sheets of the surface approach the equatorial plane, \(X_3^\star\) carries less and less of the weight, and in the limit of the doubly covered disk, a map into \(\mathbb B^2\), \(\Lambda\) has rank two. Strict complementarity thus degenerates along the family, in step with the geometry, and for \(b=8\) the certificate is already nearly singular.

The lower two blocks of Table~\ref{tab:sanity-checks} report the residuals of Section~\ref{ssec:verification}, boundary quantities as maxima over the \(bN_\ell\) quadrature nodes and interior quantities over the grid of Section~\ref{ssec:visualization}. The sphere residual \(\mathcal R_{\mathrm{sph}}\) is below \(2\times10^{-14}\) for every \(b\): with \(N_\rho\) Fourier modes per component, stationarity in the weight enforces \(|X^\star|=1\) only against trigonometric polynomials of degree \(\leq N_\rho\), and at the orders of Table~\ref{tab:sanity-checks} the tail of the Fourier series of \(|X^\star|^2-1\) is below rounding. The conformality defect and the mean curvature are at the level of the eigensolver error, \(\max_{\partial\Omega^\star}\mathcal C_{X^\star}\leq10^{-11}\) for every \(b\) and \(\max|H|\leq7\times10^{-12}\) for \(b\leq5\) and \(5\times10^{-9}\) for \(b=8\), while \(\min J\) falls from \(2.2\) for \(b=1\) to \(1.5\times10^{-4}\) for \(b=8\) as the necks shrink. The boundary length \(L_\partial\) in \eqref{eq:boundary-length} exceeds \(\bar\sigma_1^\star\) by less than \(10^{-12}\) for every \(b\). In sum, each \(X^\star\) is an immersion with \(\min J>0\), satisfies the free boundary conditions to better than \(10^{-11}\), and is minimal to better than \(5\times10^{-9}\).

\subsection{Collinear free boundary minimal tori}
\label{ssec:surfaces}

The collinear class falls within the equivariant theory of \textcite{karpukhin2024}. With \(\iota\) the reflection in the geodesic through the hole centers, each \((\Omega^\star,\langle\iota\rangle)\) is a basic reflection surface in their sense \cite[Definition~5.4]{karpukhin2024}: the quotient by \(\iota\) is an annulus, and the hole boundaries together with the fixed arcs between them form one of its boundary curves. For such pairs, \cite[Lemma~5.19]{karpukhin2024} gives \(\bar\sigma_1<4\pi\) for every \(\iota\)-invariant metric, the bound that the collinear values of Table~\ref{tab:sanity-checks} approach from below and that the maximizers of Section~\ref{ssec:equilateral}, none of which lies in the class, exceed; \cite[Propositions~5.12 and~5.15]{karpukhin2024} bound the multiplicity of \(\sigma_1\) by three, the size of the clusters of Table~\ref{tab:sanity-checks}; and \cite[Lemma~5.37]{karpukhin2024} states that a free boundary minimal immersion by first eigenfunctions with this symmetry is an embedding into \(\mathbb B^3\), two-sheeted over the equatorial disk away from the fixed curves of \(\iota\), which are convex curves in the equatorial plane, with area below \(2\pi\): this is the shape of every surface in Figure~\ref{fig:surfaces}. Their Theorem~9.15(6) gives, for every \(b\), an \(\iota\)-invariant maximizer realized by such a surface. The collinear configurations, maximizers within their symmetry class (Section~\ref{ssec:equilateral}), are the candidates for these maximizers, so that the question raised in \cite[\S1.6]{karpukhin2024} for genus zero, whether the surfaces obtained by equivariant maximization coincide with those of \textcite{Folha_Pacard_Zolotareva_2017}, receives for the tori a numerical answer in the affirmative within this class.

Figure~\ref{fig:surfaces} shows \(X^\star(\overline{\Omega^\star})\) for each \(b\), rotated so that the symmetry axes coincide with the coordinate axes (Section~\ref{ssec:visualization}). For \(b=1\) the surface is embedded, has connected boundary, and has genus one, so it is a surface of the type whose existence was proved by \textcite{Carlotto_2022}. Their surface \(M_1\) has dihedral symmetry \(\mathbb D_2\) and contains two horizontal straight lines, the fixed sets of the rotations by \(\pi\) about the horizontal axes; in \(\Omega^\star\) these are the two diagonal geodesics through the center of the hole, whose reflections have a single fixed curve, while the horizontal and vertical geodesics, fixed by basic reflections in the sense of \cite{karpukhin2024}, map to planar convex curves. The computed surface has the symmetry group of order eight noted in Section~\ref{ssec:extremal-configurations}, which contains \(\mathbb D_2\), and \(X^\star\) maps the two diagonal geodesics to straight segments to within \(2\times10^{-15}\) of the radius of the ball, perpendicular to each other and in a common horizontal plane, while the horizontal and vertical fixed curves are planar to \(10^{-9}\) but at distance \(0.1\) and \(0.3\) from the best-fitting line, so we identify it with \(M_1\).

For \(b\geq2\) the surfaces consist of two sheets, each close to a flat disk, joined by a catenoidal neck at the center and by \(b\) necks equally spaced along the equator; the sheets approach each other as \(b\) grows and for \(b=8\) are nearly coincident. This is the shape of the surfaces constructed by \textcite{Folha_Pacard_Zolotareva_2017} for \(b\) large, whose symmetry group is generated by the rotation by \(2\pi/b\) about the vertical axis, the reflection in the horizontal plane, which exchanges the two sheets, and the reflection in a vertical plane through a neck; for \(b=2\) it is the genus one critical catenoid of \textcite{Franz_Ketover_Schulz_2024}, with the same symmetry group of order eight, whose first Steklov eigenvalue they prove to be one. The computed surfaces have these symmetries to rounding, \(2\times10^{-13}\) for \(\omega^\star\) and for \(X^\star\) on \(\partial\Omega^\star\) under the generators, as a symmetric initialization and the invariance of the iteration ensure; in \(\Omega^\star\) the rotation is the translation by \(1/b\) along the geodesic, the horizontal reflection is the reflection in the geodesic, and the vertical reflection is the reflection in a line through a hole center perpendicular to the geodesic.

\subsection{Maximizers and saddle points}
\label{ssec:equilateral}

In this subsection the symmetry restriction of Section~\ref{ssec:setup} is not imposed; the optimization is initialized at the collinear configurations, perturbed in symmetry-breaking directions, and at random data; Table~\ref{tab:equilateral} lists the configurations with the largest values found.
\begin{figure}[p]
\centering
\includegraphics[width=\linewidth,height=0.84\textheight,keepaspectratio]{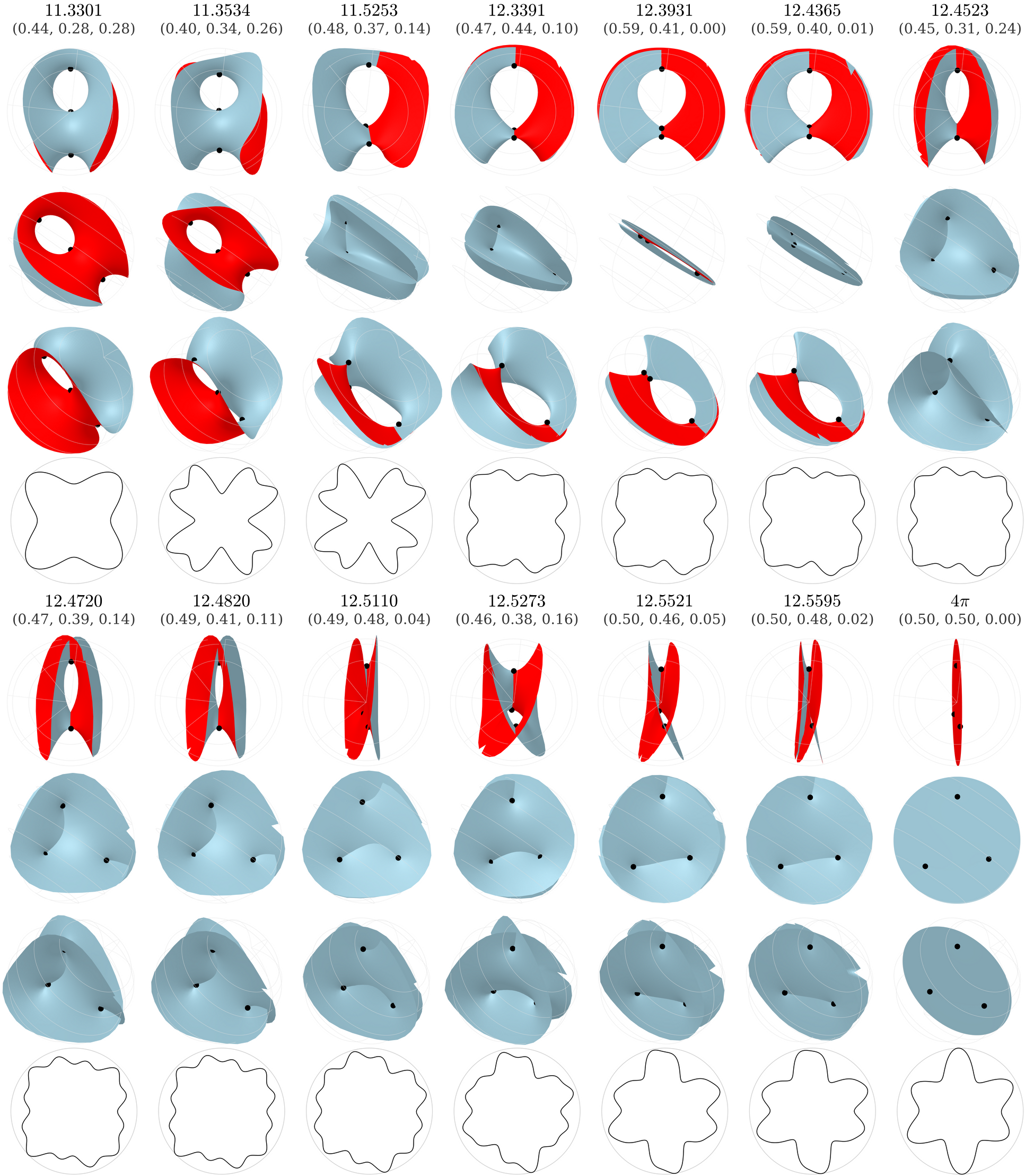}
\caption{Fourteen iterates of the unrestricted optimization initialized at the \(b=1\) collinear configuration of Table~\ref{tab:sanity-checks}, perturbed by \(10^{-3}\) in a symmetry-breaking direction, ending at the maximizer, a doubly covered disk (Section~\ref{ssec:equilateral}); \(\bar\sigma_1\) and the eigenvalues of the multiplier \(\Lambda\) are given above each column. From top to bottom: the spectral map in three views with a fixed camera and scale, and the weight \(\rho/\max\rho\). Red and gray distinguish the two sides of the surface. Black dots mark the images of the three half-periods relative to the hole center, which are the ramification points at the maximizer.}
\label{fig:b1-fixed-camera}
\end{figure}

\begin{table}[t]
\centering
\footnotesize
\setlength{\tabcolsep}{4pt}
\renewcommand{\arraystretch}{0.95}

\sisetup{round-mode=places, round-precision=1, table-format=-1.1e+2, group-minimum-digits=20}
\begin{tabular}{@{}l *{6}{S} @{}}
\toprule
& {\(b=1\)} & {\(b=2\)} & {\(b=3\)} & {\(b=4\)} & {\(b=5\)} & {\(b=8\)}\\
\midrule
\multicolumn{7}{@{}l}{\textit{Computed configuration}}\\
\addlinespace[1pt]
\(N_\rho\) & {56} & {48} & {56} & {64} & {80} & {120}\\
\(m\) & {2} & {5} & {6} & {6} & {6} & {6}\\
\(\operatorname{rank}\Lambda^\star\) & {2} & {5} & {6} & {6} & {6} & {6}\\
\(\bar\sigma_1^\star\) & {\tablenum[scientific-notation=false,round-mode=none,table-format=2.10]{12.5663706144}} & {\tablenum[scientific-notation=false,round-mode=none,table-format=2.10]{15.4719507709}} & {\tablenum[scientific-notation=false,round-mode=none,table-format=2.10]{18.4863452307}} & {\tablenum[scientific-notation=false,round-mode=none,table-format=2.10]{20.4823390210}} & {\tablenum[scientific-notation=false,round-mode=none,table-format=2.10]{21.8816924184}} & {\tablenum[scientific-notation=false,round-mode=none,table-format=2.10]{25.9226908749}} \\
\(\Im\tau^\star\) & {\tablenum[scientific-notation=false,round-mode=none,table-format=2.10]{0.8730328765}} & {\tablenum[scientific-notation=false,round-mode=none,table-format=2.10]{0.8660254038}} & {\tablenum[scientific-notation=false,round-mode=none,table-format=2.10]{0.8660254038}} & {\tablenum[scientific-notation=false,round-mode=none,table-format=2.10]{0.8660254038}} & {\tablenum[scientific-notation=false,round-mode=none,table-format=2.10]{0.8660254038}} & {\tablenum[scientific-notation=false,round-mode=none,table-format=2.10]{0.8660530352}}\\
\(r^\star\) & {\tablenum[scientific-notation=false,round-mode=none,table-format=2.10]{0.3100000000}} & {\begin{tabular}[t]{@{}c@{}}\tablenum[scientific-notation=false,round-mode=none,table-format=2.10]{0.2926597784} \\ \tablenum[scientific-notation=false,round-mode=none,table-format=2.10]{0.0583025969}\end{tabular}} & {\begin{tabular}[t]{@{}c@{}}\tablenum[scientific-notation=false,round-mode=none,table-format=2.10]{0.1180160745} \\ \((\times3)\)\end{tabular}} & {\begin{tabular}[t]{@{}c@{}}\tablenum[scientific-notation=false,round-mode=none,table-format=2.10]{0.0928990944} \\ \((\times4)\)\end{tabular}} & {\begin{tabular}[t]{@{}c@{}}\tablenum[scientific-notation=false,round-mode=none,table-format=2.10]{0.0767593548} \\ \((\times3)\) \\ \tablenum[scientific-notation=false,round-mode=none,table-format=2.10]{0.0850527093} \\ \tablenum[scientific-notation=false,round-mode=none,table-format=2.10]{0.1129128908}\end{tabular}} & {\begin{tabular}[t]{@{}c@{}}\tablenum[scientific-notation=false,round-mode=none,table-format=2.10]{0.0635202457} \\ \((\times8)\)\end{tabular}} \\
\midrule
\multicolumn{7}{@{}l}{\textit{Stationarity (Section~\ref{ssec:two-stage-optimization})}}\\
\addlinespace[1pt]
\(\delta_m(\mathbf p^\star)\) \eqref{eq:cluster-spread} & 1.48849959e-13 & 1.15959547e-14 & 9.99338209e-15 & 3.36497923e-14 & 3.11731384e-14 & 3.48107872e-14 \\
\(\mathcal K(\mathbf p^\star,\Lambda^\star)\) \eqref{eq:kkt-residual} & 6.92557123e-13 & 4.48252546e-14 & 6.18219777e-13 & 2.90253932e-14 & 2.43138842e-14 & 2.05280237e-13 \\
\(\lambda_{\min}(\Lambda^\star)\) & 5.0000000000e-01 & 9.2705821483e-02 & 5.2370983753e-02 & 1.1713774232e-01 & 6.4092046144e-02 & 1.6111811386e-01\\
\(\bar\sigma_{m+1}/\bar\sigma_1\) & {\tablenum[scientific-notation=false,round-mode=none,table-format=2.10]{1.0154134605}} & {\tablenum[scientific-notation=false,round-mode=none,table-format=2.10]{1.7498710820}} & {\tablenum[scientific-notation=false,round-mode=none,table-format=2.10]{1.3839616754}} & {\tablenum[scientific-notation=false,round-mode=none,table-format=2.10]{1.5046259356}} & {\tablenum[scientific-notation=false,round-mode=none,table-format=2.10]{1.5669135947}} & {\tablenum[scientific-notation=false,round-mode=none,table-format=2.10]{1.8530070343}} \\
\midrule
\multicolumn{7}{@{}l}{\textit{Free boundary conditions on \(\partial\Omega^\star\) (Section~\ref{ssec:verification})}}\\
\addlinespace[1pt]
\(\max\mathcal R_{\mathrm{sph}}\) & 4.79616347e-14 & 6.23945340e-14 & 2.70894418e-13 & 4.21884749e-15 & 7.10542736e-15 & 3.13415960e-13 \\
\(L_\partial-\bar\sigma_1^\star\) \eqref{eq:boundary-length} & 9.89430760e-13 & 8.34887715e-14 & 1.10134124e-13 & 2.84217094e-13 & 3.33955086e-13 & 3.41060513e-13 \\
\(\max_{\partial\Omega^\star}\mathcal C_{X^\star}\) \eqref{eq:conformality-defect} & 3.35238867e-13 & 1.97082232e-12 & 4.41624863e-12 & 4.63393464e-13 & 4.24653195e-13 & 4.87663724e-12 \\
\midrule
\multicolumn{7}{@{}l}{\textit{Minimality and immersion in \(\Omega^\star\) (Section~\ref{ssec:verification})}}\\
\addlinespace[1pt]
\(\max|H|\) \eqref{eq:mean-curvature} & {--} & 1.87465711e-12 & 2.06076529e-12 & 3.91436677e-14 & 1.49493695e-13 & 3.11399663e-12 \\
\(\min J\) & {--} & 3.3550432275e+00 & 4.6744131702e+00 & 5.4184166790e+00 & 6.2958815851e+00 & 1.0458619549e+01\\
\bottomrule
\end{tabular}
\caption{The configurations with the largest values of \(\bar\sigma_1\) found (Section~\ref{ssec:equilateral}), computed at \(N_\ell=400\) with the listed \(N_\rho\): cluster size, rank of the certificate, extremal value, imaginary part of the reduced modulus (\(\Re\tau^\star=\tfrac12\) for every \(b\)), radii, stationarity measures, gap to the next eigenvalue, and residuals in their \(m\)-dimensional form (conformality defect from the singular values of \(dX^\star\), length of the mean curvature vector). For \(b=1\) the spectral map is a branched double cover of the disk, with \(J=0\) at the three ramification points, where \(\mathcal C_{X^\star}\) and \(H\) are undefined.}
\label{tab:equilateral}
\end{table}

\emph{Saddle points.} At each collinear configuration we computed the Hessian of the Lagrangian of \eqref{eq:local-sqp} on the tangent space of the coalescence manifold, by finite differences of the derivatives \eqref{eq:cluster-derivative}. It has one positive eigenvalue for \(b=1\), three for \(b=2\) and five for \(b=3\), all in symmetry-breaking directions: for \(b=1\) the modulus leaves the imaginary axis together with an odd Fourier mode of the weight, for \(b=2\) the two holes move off the geodesic in opposite directions and their radii and weights become unequal, for \(b=3\) the centers leave the geodesic; every direction that preserves the symmetry has negative curvature. A straight perturbation in such a direction lowers \(\bar\sigma_1\), because it splits the triple eigenvalue. Accordingly, a perturbation of size \(10^{-3}\) followed by the iteration of Section~\ref{ssec:optimizer-settings}, which re-coalesces the cluster at every step, raises \(\bar\sigma_1\) monotonically from the collinear values to \(4\pi\), \(18.49\), \(20.48\) and \(21.88\) for \(b=1,3,4,5\), and for \(b=2\) to \(15.37\), a local maximizer on the square torus with a four-fold cluster, below the value of Table~\ref{tab:equilateral}. The collinear configurations are therefore maximizers within their symmetry class only, and it is the restriction of Section~\ref{ssec:setup} that makes them reachable; for \(b=1\) this answers the question left open in \cite{Carlotto_2022}.

\emph{The case \(b=1\).} The Fraser--Schoen bound \(\bar\sigma_1\leq2\pi(\gamma+b)=4\pi\) \cite[Theorem~2.2]{Fraser_2015} is attained: the configuration of Table~\ref{tab:equilateral} has \(\bar\sigma_1^\star=4\pi\) to ten digits, and Figure~\ref{fig:b1-fixed-camera} shows how the perturbed collinear configuration reaches it. The surface first flattens: at \(\bar\sigma_1=12.39\) the multiplier has rank two and the map is planar, folded onto itself, with the images of the three half-periods \(c+\tfrac12\), \(c+\tfrac\tau2\), \(c+\tfrac{1+\tau}2\) relative to the hole center on a line. The smallest eigenvalue of \(\Lambda\) then returns to \(0.24\) at \(\bar\sigma_1=12.45\) before vanishing again at the maximizer, the images of the half-periods spread from a line to a nearly equilateral triangle, and the weight changes from a four-fold to a six-fold profile. The maximizer is not a surface in \(\mathbb B^3\): two eigenvalues coincide with \(\bar\sigma_1\) and the third lies \(1.5\%\) above, \(\Lambda^\star\) has rank two with eigenvalues \(\tfrac12\), \(\tfrac12\), and the spectral map \(X^\star=(u_1,u_2)\) is a branched double cover of the unit disk by the torus. It maps the boundary twice around the unit circle (\(\max||X^\star|-1|\leq10^{-13}\), Table~\ref{tab:equilateral}; boundary length \(4\pi=\bar\sigma_1^\star\)), and \(dX^\star\) vanishes exactly at the three half-periods, the fixed points of the half-turn of the torus about the hole center, which is the deck transformation of the cover (Figure~\ref{fig:b1-fixed-camera}); three ramification points are what the Riemann--Hurwitz formula requires. This is the equality case of the bound of \cite{Fraser_2015}: \(X^\star\) is the holomorphic map of degree \(\gamma+b=2\) onto the disk from their proof, with first eigenfunctions as coordinates, and \(\sigma^*(1,1)=4\pi\). The maximum is attained on a family rather than at one configuration: the maximizers reached from different initializations have \(|\Re\tau^\star|\in[0.39,0.50]\), \(|\tau^\star|\in[1.00,1.04]\) and \(r^\star\in[0.31,0.33]\), all with \(\bar\sigma_1=4\pi\); the moduli lie near the hexagonal modulus \(e^{i\pi/3}\) and never on the imaginary axis, in accordance with the bound \(\bar\sigma_1<4\pi\) of \cite[Lemma~5.19]{karpukhin2024} for metrics invariant under a reflection whose quotient is an annulus (Section~\ref{ssec:surfaces}).

\emph{The cases \(b\geq2\).} The largest values found are at configurations of a different kind: for \(b=2\) a five-fold cluster on the equilateral torus with holes of radii \(0.293\) and \(0.058\); for \(b=3\) a six-fold cluster, equilateral, with three equal holes equispaced on a closed geodesic; for \(b=4\) six-fold, equilateral, four equal holes in two rows; for \(b=5\) six-fold, equilateral, unequal holes; for \(b=8\) six-fold, eight equal holes. Their certificates have full rank \(m\) and their residuals in the \(m\)-dimensional form of Section~\ref{ssec:verification} are listed in Table~\ref{tab:equilateral}, so by Theorem~\ref{prop:certificate} each spectral map is a free boundary minimal immersion of the torus into \(\mathbb B^5\) or \(\mathbb B^6\), none a surface in \(\mathbb B^3\); for \(b=2,3\) the Hessian on the tangent space is negative definite, so these two are local maximizers. For \(b=2,\dots,5\) they lie above the lower bounds \(\sigma^*(1,b)>\sigma^*(0,b+1)=12.01\), \(13.67\), \(14.47\), \(15.43\) that follow from the genus-zero values of \textcite[Table~2]{Kao_2020} and the gluing inequality of \cite[Theorem~1.3]{Matthiesen_Petrides_2020}, and below the Fraser--Schoen bound \(2\pi(b+1)\); whether they attain \(\sigma^*(1,b)\) we cannot decide.

As \(b\to\infty\), \textcite{Karpukhin_Stern_2021} show that the maximal value of \(\bar\sigma_1\) on a surface of genus \(\gamma\) converges to the maximal area-normalized first Laplace eigenvalue of the closed surface of genus \(\gamma\), and that the corresponding free boundary minimal surfaces converge to a closed minimal surface in the sphere; for the torus the limit is governed by the equilateral torus \cite{Nadirashvili_1996}, whose first eigenspace is six-dimensional. The family of \textcite{Folha_Pacard_Zolotareva_2017} computed here converges instead to the doubly covered disk (Section~\ref{ssec:extremal-configurations}), and \(\Im\tau^\star\) increases with \(b\) (Figure~\ref{fig:eigenfunction_contours}), away from the equilateral modulus \(e^{i\pi/3}\) rather than toward it; the maximizers of Table~\ref{tab:equilateral} lie outside this family for every \(b\), and the six-fold clusters of \(b\geq3\) already match the dimension of the first eigenspace of the equilateral torus that governs the limit of \textcite{Karpukhin_Stern_2021}.

\section{Discussion}
\label{sec:discussion}

We have computed free boundary minimal tori by maximizing the first normalized Steklov eigenvalue over the modulus of a flat torus, the centers and radii of \(b\) round holes, and a boundary weight, a parametrization that covers the moduli space of conformal structures. The weighted Steklov problem is a generalized eigenvalue problem for a pencil of boundary integral operators whose left eigenvector is the trace of the eigenfunction (Theorem~\ref{thm:weighted-steklov-bie}), so that the derivatives with respect to all parameters, the modulus included, follow from one formula (Proposition~\ref{prop:unified-derivative}), and the coalescence of eigenvalues at a maximizer is handled by a two-stage nonsmooth method (Section~\ref{ssec:two-stage-optimization}) whose stationarity certificate rescales the eigenfunctions into the free boundary minimal immersion (Theorem~\ref{prop:certificate}). The surfaces of \textcite{Carlotto_2022} and of \textcite{Folha_Pacard_Zolotareva_2017} are critical points of \(\bar\sigma_1\), maximizers within their symmetry class but saddle points; for connected boundary the supremum is the Fraser--Schoen bound \(4\pi\), attained by branched double covers of the disk, and for \(b\geq2\) the largest values found are attained by free boundary minimal tori in \(\mathbb B^5\) and \(\mathbb B^6\). Beyond genus zero, then, the free boundary minimal surfaces in \(\mathbb B^3\) produced by the Fraser--Schoen correspondence are critical points rather than maximizers, at least for these topological types.

The method does not require the symmetry class of Section~\ref{ssec:setup}. Prescribing a symmetry group splits the eigenproblem into sectors with simple eigenvalues and turns \eqref{eq:intro-three-dimensional-cluster} into the equality of sector eigenvalues, the smooth system solved in the refinement step of Section~\ref{ssec:optimizer-settings}; in the equivariant setting of \cite{karpukhin2024} this is the natural route to lower-symmetry configurations and two-dimensional arrangements of holes. The method also applies to higher eigenvalues \(\bar\sigma_k\), \(k\geq2\) \cite{Fraser_2020,Kao_2020}.

The residuals of Section~\ref{sec:experiments} are numerical evidence, not proof. It would be interesting to prove convergence of the method, or to turn a computed configuration into a computer-assisted existence proof; for the connected-boundary torus this would complement the variational proof of \cite{Carlotto_2022}. Our criteria certify immersions, and embeddedness is judged from the plots. For \(b=1\) the maximizer is not even an immersion (Section~\ref{ssec:equilateral}), although the class contains an embedded free boundary minimal surface, so already for the torus with connected boundary the Fraser--Schoen correspondence produces a branched map at the maximum; which classes, or which equivariant subclasses \cite{karpukhin2024}, have embedded maximizers is the question this raises.

\subsection*{Data availability.}

The code and parameters of the computed configurations are available in the repository \cite{Chen_Osting_2026_data}.

\subsection*{Acknowledgments}

B.O. acknowledges useful discussions with Chiu-Yen Kao, Micka\"el Nahon, and \'Edouard Oudet.
The authors used Claude (Anthropic) to assist with editing the text, verifying references, and checking derivations. ChatGPT (OpenAI) was used in developing parts of the numerical code. All AI-assisted text and code were checked by the authors, who take full responsibility for the content of the paper.

\printbibliography
\clearpage
\appendix

\section*{Appendices}

\addcontentsline{toc}{section}{Appendices}

\section{The modulus derivative of the torus Green function}
\label{sec:analytic formula for dtau}

We derive the Wirtinger derivative \eqref{eq:dtau-G-closed} of the torus Green function \eqref{eq:torus-green-function}, show that it is smooth, and record the series by which it is evaluated.

\smallskip\noindent\emph{The logarithmic term.} Fix \(z\notin\mathcal L_\tau\). The theta function \(\vartheta_1(z;\tau)\) is holomorphic in \(\tau\in\mathbb H\), so \(\partial_{\bar\tau}\vartheta_1=0\) and, taking complex conjugates, \(\partial_\tau\overline{\vartheta_1}=0\). Writing \(\log|\vartheta_1|=\tfrac12\bigl(\log\vartheta_1+\overline{\log\vartheta_1}\bigr)\) and applying \(\partial_\tau\), the conjugate term drops out:
\begin{equation}\label{eq:AppA_wirtinger_logabs}
\partial_\tau\log|\vartheta_1(z;\tau)|=\tfrac12\,\partial_\tau\log\vartheta_1(z;\tau)=\frac{\partial_\tau\vartheta_1(z;\tau)}{2\,\vartheta_1(z;\tau)}.
\end{equation}
This is where the Wirtinger convention matters: \(G_\tau\) is real, so \(\partial_{\bar\tau}G_\tau=\overline{\partial_\tau G_\tau}\), and the real derivatives are recovered from \(\partial_\tau G_\tau\) alone as \(\partial_{\tau_1}G_\tau=2\Re\,\partial_\tau G_\tau\) and \(\partial_{\tau_2}G_\tau=-2\Im\,\partial_\tau G_\tau\).

\smallskip\noindent\emph{The quadratic term.} Since \(\Im\tau=(\tau-\bar\tau)/(2i)\), we have \(\partial_\tau\Im\tau=\tfrac1{2i}\), hence \(\partial_\tau\bigl[(\Im z)^2/(2\Im\tau)\bigr]=-\tfrac{(\Im z)^2}{2(\Im\tau)^2}\cdot\tfrac1{2i}=\tfrac{i(\Im z)^2}{4(\Im\tau)^2}\). Combining with \eqref{eq:AppA_wirtinger_logabs} and the factor \(-\tfrac1{2\pi}\) in \eqref{eq:torus-green-function} gives \eqref{eq:dtau-G-closed}.

\smallskip\noindent\emph{Smoothness at the lattice points.} In the normalization of \eqref{eq:torus-green-function}, in which the lattice is \(\mathbb Z+\tau\mathbb Z\), the theta function satisfies the heat equation
\begin{equation}\label{eq:AppA_heat}
\partial_\tau\vartheta_1(z;\tau)=-\frac{i}{4\pi}\,\partial_z^2\vartheta_1(z;\tau)
\end{equation}
\cite[\S20.13]{DLMF}. Since \(\vartheta_1\) is odd in \(z\) with simple zeros exactly at the lattice points, \(\partial_z^2\vartheta_1\) is odd as well and vanishes there, so the quotient \(\partial_\tau\vartheta_1/\vartheta_1\) in \eqref{eq:dtau-G-closed} has removable singularities at \(\mathcal L_\tau\), with limit \(\partial_{\tau z}\vartheta_1/\partial_z\vartheta_1\) by l'H\^opital's rule. Consequently \(\partial_\tau G_\tau\) extends to a smooth function on all of \(\mathbb C\): the logarithmic singularity \(-\tfrac1{2\pi}\log|z|\) of \(G_\tau\) is independent of \(\tau\) and is annihilated by \(\partial_\tau\). This is why the operators \(\dot{\mathcal S}\) and \(\dot{\mathcal D}^*\) of Proposition~\ref{prop:unified-derivative} require no singular quadrature (Section~\ref{ssec:nystrom}).

\smallskip\noindent\emph{Series for the implementation.} With \(q:=e^{\pi i\tau}\) and \(\vartheta_1(z;\tau)=2\sum_{n=0}^\infty(-1)^nq^{(n+1/2)^2}\sin((2n+1)\pi z)\), only the factor \(q^{(n+1/2)^2}\) depends on \(\tau\), and \(\partial_\tau q^{(n+1/2)^2}=(n+\tfrac12)^2\,\pi i\,q^{(n+1/2)^2}\), so term-by-term differentiation gives
\begin{equation}\label{eq:K_series_impl}
\partial_\tau\vartheta_1(z;\tau)=2\pi i\sum_{n=0}^\infty(-1)^n\bigl(n+\tfrac12\bigr)^2q^{(n+1/2)^2}\sin\bigl((2n+1)\pi z\bigr);
\end{equation}
the series for \(\partial_z\vartheta_1\) and \(\partial_{\tau z}\vartheta_1\), needed for the diagonal limit, are obtained in the same way. All four series converge geometrically, since \(|q|=e^{-\pi\Im\tau}<1\), and are summed together to a prescribed relative tolerance (Section~\ref{ssec:nystrom}).

\smallskip\noindent\emph{Expression through Weierstrass functions.} Dividing \eqref{eq:AppA_heat} by \(\vartheta_1\) gives \(\partial_\tau\log\vartheta_1=-\tfrac{i}{4\pi}\bigl(\partial_z^2\log\vartheta_1+(\partial_z\log\vartheta_1)^2\bigr)\). Two classical identities express the right-hand side through the Weierstrass functions \(\wp\), \(\zeta_{\mathrm W}\) of the lattice \(\mathcal L_\tau\), the quasi-period \(\eta_1\) and the Eisenstein series \(E_2(\tau)=1-24\sum_{n=1}^\infty nq^{2n}/(1-q^{2n})\): \(\partial_z\log\vartheta_1=\zeta_{\mathrm W}(z)-2\eta_1z\), which follows from \(\sigma_{\mathrm W}(z)=C(\tau)\,e^{\eta_1z^2}\vartheta_1(z;\tau)\) \cite[\S23.6]{DLMF}, and \(-\partial_z^2\log\vartheta_1=\wp(z)+\tfrac{\pi^2}{3}E_2(\tau)\), obtained by comparing the Laurent expansions at \(z=0\) of \(\partial_z^2\log\vartheta_1\) \cite[Eq.~20.5.10]{DLMF} and of \(\wp\). Substituting into \eqref{eq:dtau-G-closed},
\begin{equation}\label{eq:AppA_final_dtau_G}
\partial_\tau G_\tau(z)=-\frac{i}{16\pi^2}\wp(z)+\frac{i}{16\pi^2}\bigl(\zeta_{\mathrm W}(z)-2\eta_1z\bigr)^2-\frac{i}{48}E_2(\tau)+\frac{i(\Im z)^2}{4(\Im\tau)^2}.
\end{equation}
The same identities give the Hessian of \(G_\tau\) used in Section~\ref{ssec:analytic-sensitivities}: for the harmonic function \(\log|\vartheta_1|=\Re\log\vartheta_1\) one has \(\partial_{xx}=-\partial_{yy}=\Re f''\) and \(\partial_{xy}=-\Im f''\) with \(f'':=\partial_z^2\log\vartheta_1=-\wp(z)-\tfrac{\pi^2}{3}E_2(\tau)\), and the quadratic term of \eqref{eq:torus-green-function} adds \(1/\Im\tau\) to \(\partial_{yy}G_\tau\). Finally, the constant \(-\tfrac{i}{48}E_2(\tau)\) in \eqref{eq:AppA_final_dtau_G} integrates to zero against the zero-mean densities on which \(\dot{\mathcal S}\) and \(\dot{\mathcal D}^*\) act.

\section{Close evaluation by quadrature by expansion}
\label{sec:qbx}

When a target \(z\in\Omega^\star\) lies within a few node spacings of a source circle \(\Gamma_\ell\), the trapezoidal rule applied to \eqref{eq:eigenfunction-layer-representation} loses accuracy, because the integrand varies on the scale of the distance to the boundary. We then evaluate the contribution of \(\Gamma_\ell\) by quadrature by expansion (QBX) \cite{Klockner2013}, in the standard form for the logarithmic kernel; this appendix records the formulas and parameters used.

\smallskip\noindent\emph{Kernel split.} With the period-one normalization of \(\vartheta_1\), the torus Green function splits as
\begin{equation}\label{eq:sqbx-kernel-split}
G_\tau(w)=-\frac{1}{2\pi}\log|w|+R_\tau(w),\qquad R_\tau(w):=-\frac{1}{2\pi}\log\Bigl|\frac{\vartheta_1(w;\tau)}{w}\Bigr|+\frac{(\Im w)^2}{2\Im\tau},
\end{equation}
where \(R_\tau\) is smooth near \(w=0\), since \(\vartheta_1\) has a simple zero there (Appendix~\ref{sec:analytic formula for dtau}); its singularities at the other lattice points are not removed, which is why the construction is applied to one periodic copy of one component at a time. By \eqref{eq:modified-single-layer-potential}--\eqref{eq:eigenfunction-layer-representation}, \(u=\overline\varphi+\sum_{\ell\in[b]}U_\ell\) with \(\overline\varphi:=M[\varphi]\), and each component potential splits accordingly into \(U_\ell=U_\ell^{\mathrm{sing}}+U_\ell^{\mathrm{reg}}\),
\[U_\ell^{\mathrm{sing}}(z):=-\frac{1}{2\pi}\int_{\Gamma_\ell}(\varphi_\ell-\overline\varphi)\log|z-\xi|\,ds_\xi,\qquad
U_\ell^{\mathrm{reg}}(z):=\int_{\Gamma_\ell}(\varphi_\ell-\overline\varphi)\,R_\tau(z-\xi)\,ds_\xi .\]
The regular part is integrated by the trapezoidal rule; QBX is applied to \(U_\ell^{\mathrm{sing}}\) only.

\smallskip\noindent\emph{Expansion.} Let \(\xi_0\in\Gamma_\ell\) be the boundary point nearest to \(z\), so that \(z=\xi_0-|z-\xi_0|\,\nu_\tau(\xi_0)\), the normal \(\nu_\tau\) pointing into the removed disk. Choose \(\delta_\ell>|z-\xi_0|/2\) and set \(\zeta:=\xi_0-\delta_\ell\,\nu_\tau(\xi_0)\), the center of the disk of radius \(\delta_\ell\) tangent to \(\Gamma_\ell\) at \(\xi_0\) from the exterior side; \(\delta_\ell\) is taken small enough that this disk avoids all other periodic copies of \(\Gamma_\ell\). We take \(\delta_\ell=\max\{\min(4|z-\xi_0|,\,2h_\ell),\,0.75\,|z-\xi_0|\}\), with \(h_\ell=2\pi r_\ell/N_\ell\) the node spacing on \(\Gamma_\ell\). Writing \(w_\ell:=(z-\zeta)/\delta_\ell\), we have \(|w_\ell|<1\) and \(|\xi-\zeta|\geq\delta_\ell\) for all \(\xi\in\Gamma_\ell\), so \(\log|z-\xi|=\log|\xi-\zeta|+\Re\log\bigl(1-(z-\zeta)/(\xi-\zeta)\bigr)\) can be expanded in a geometric series that converges uniformly in \(\xi\). Integrating term by term,
\begin{equation}\label{eq:sqbx-exact-expansion}
U_\ell^{\mathrm{sing}}(z)=A_{\ell,0}+\sum_{n=1}^\infty\bigl(A_{\ell,n}w_\ell^{\,n}+B_{\ell,n}\overline{w_\ell}^{\,n}\bigr),
\end{equation}
with coefficients
\[A_{\ell,0}=-\frac{1}{2\pi}\int_{\Gamma_\ell}(\varphi_\ell-\overline\varphi)\log|\xi-\zeta|\,ds_\xi,\quad
A_{\ell,n}=\frac{1}{4\pi n}\int_{\Gamma_\ell}(\varphi_\ell-\overline\varphi)\Bigl(\frac{\delta_\ell}{\xi-\zeta}\Bigr)^{n}ds_\xi,\quad
B_{\ell,n}=\overline{A_{\ell,n}},\]
for \(n\geq1\), the last identity because \(\varphi_\ell\) is real. The truncation error after \(Q\) terms is bounded by \(\|\varphi_\ell-\overline\varphi\|_{L^1(\Gamma_\ell)}\,|w_\ell|^{Q+1}/\bigl(2\pi(Q+1)(1-|w_\ell|)\bigr)\), so \eqref{eq:sqbx-exact-expansion} converges geometrically in \(Q\).

\smallskip\noindent\emph{Discretization.} The coefficient integrals involve only smooth integrands, but they vary on the scale \(\delta_\ell\); we therefore interpolate the computed nodal density onto a finer uniform grid of \(N_\ell^{\mathrm Q}=64N_\ell\) nodes on \(\Gamma_\ell\) by Fourier interpolation and apply the trapezoidal rule \cite{Barnett_2014}, which does not alter the underlying boundary integral solution. With \(Q=12\) terms, the contribution of a close component is
\[U_\ell(z)\approx \widehat A_{\ell,0}+\sum_{n=1}^Q\bigl(\widehat A_{\ell,n}w_\ell^{\,n}+\overline{\widehat A_{\ell,n}}\,\overline{w_\ell}^{\,n}\bigr)+h_\ell^{\mathrm Q}\sum_i\bigl[\varphi_\ell(\xi_i)-\overline\varphi\bigr]R_\tau(z-\xi_i),\]
where \(\widehat A_{\ell,n}\) are the trapezoidal approximations of the coefficients and \(h_\ell^{\mathrm Q}:=2\pi r_\ell/N_\ell^{\mathrm Q}\) the spacing of the fine grid. The gradient and Hessian of \(u\) are obtained by differentiating the expansion in \(z\), \(\partial_zU_\ell^{\mathrm{sing}}\approx\delta_\ell^{-1}\sum_{n=1}^Q n\,\widehat A_{\ell,n}w_\ell^{\,n-1}\) and \(\partial_z^2U_\ell^{\mathrm{sing}}\approx\delta_\ell^{-2}\sum_{n=2}^Q n(n-1)\,\widehat A_{\ell,n}w_\ell^{\,n-2}\), and by integrating \(\nabla R_\tau\) and its Hessian on the fine grid for the regular part.

\end{document}